\documentclass[a4paper,11pt, leqno]{article}
\usepackage[latin1]{inputenc}
\usepackage[francais,english]{babel}
\usepackage[intlimits]{amsmath}
\usepackage{amssymb}
\usepackage{latexsym}
\usepackage{mathrsfs}
\usepackage[all]{xy}

\RequirePackage{ifthen}
\RequirePackage{verbatim}

\renewcommand{\AA}{\mathbb{A}}
\newcommand{\BB}{\mathbb{B}}
\newcommand{\CC}{\mathbb{C}}
\newcommand{\DD}{\mathbb{D}}

\newcommand{\FF}{\mathbb{F}}
\newcommand{\GG}{\mathbb{G}}

\newcommand{\KK}{\mathbb{K}}

\newcommand{\MM}{\mathbb{M}}

\newcommand{\PP}{\mathbb{P}}
\newcommand{\QQ}{\mathbb{Q}}
\newcommand{\RR}{\mathbb{R}}
\renewcommand{\SS}{\mathbb{S}}

\newcommand{\VV}{\mathbb{V}}

\newcommand{\XX}{\mathbb{X}}

\newcommand{\ZZ}{\mathbb{Z}}

\renewcommand{\hbar}{\bar{h}}

\newcommand{\kbar}{\bar{k}}

\newcommand{\sbar}{\bar{s}}

\newcommand{\Mbar}{\bar{M}}

\newcommand{\Sbar}{\bar{S}}

\newcommand{\Ubar}{\bar{U}}

\newcommand{\hgbar}{\bar{\eta}}

\newcommand{\lgbar}{\bar{\lambda}}

\newcommand{\tgbar}{\bar{\tau}}

\newcommand{\FFbar}{\overline{\FF}}

\newcommand{\XXbar}{\overline{\XX}}

\newcommand{\SStilde}{\tilde{\SS}}

\newcommand{\Abf}{\mathbf{A}}

\newcommand{\Fbf}{\mathbf{F}}

\newcommand{\Mbf}{\mathbf{M}}
\newcommand{\Nbf}{\mathbf{N}}

\newcommand{\Vbf}{\mathbf{V}}

\newcommand{\Xbf}{\mathbf{X}}

\newcommand{\igbf}{\mathbf{\iota}}

\newcommand{\lgbf}{\mathbf{\lambda}}

\newcommand{\Acal}{{\mathcal A}}
\newcommand{\Bcal}{{\mathcal B}}

\newcommand{\Fcal}{{\mathcal F}}
\newcommand{\Gcal}{{\mathcal G}}

\newcommand{\Lcal}{{\mathcal L}}
\newcommand{\Mcal}{{\mathcal M}}
\newcommand{\Ncal}{{\mathcal N}}
\newcommand{\Ocal}{{\mathcal O}}

\newcommand{\Vcal}{{\mathcal V}}

\newcommand{\Xcal}{{\mathcal X}}

\newcommand{\Bscr}{{\mathscr B}}

\newcommand{\Fscr}{{\mathscr F}}

\newcommand{\Nscr}{{\mathscr N}}
\newcommand{\Oscr}{{\mathscr O}}

\newcommand{\dtilde}{\tilde{d}}

\newcommand{\Dtilde}{\tilde{D}}

\newcommand{\Htilde}{\tilde{H}}

\newcommand{\Jtilde}{\tilde{J}}

\newcommand{\Lgtilde}{\tilde\Lambda}

\newcommand{\eps}{\varepsilon}

\DeclareMathOperator{\BT}{BT}
\DeclareMathOperator{\Cent}{Cent}

\DeclareMathOperator{\Coker}{Coker}

\DeclareMathOperator{\End}{End}

\DeclareMathOperator{\GL}{GL}

\DeclareMathOperator{\GU}{GU}
\DeclareMathOperator{\height}{height}
\DeclareMathOperator{\Hom}{Hom}

\DeclareMathOperator{\id}{id}
\renewcommand{\Im}{\mathop{\rm Im}}

\DeclareMathOperator{\Ker}{Ker}

\DeclareMathOperator{\Lie}{Lie}

\DeclareMathOperator{\Par}{Par}

\DeclareMathOperator{\Res}{Res}
\DeclareMathOperator{\rk}{rk}
\DeclareMathOperator{\SL}{SL}

\DeclareMathOperator{\Spf}{Spf}
\DeclareMathOperator{\Spec}{Spec}
\DeclareMathOperator{\SU}{SU}

\DeclareMathOperator{\Tr}{Tr}

\newcommand{\red}{\textnormal{red}}

\newcommand{\cross}{{}^{\times}}

\newcommand{\lrangle}{\mathop{\langle\,\mathop{,}\,\rangle}}
\newcommand{\lrcurly}{\mathop{\{\,\mathop{,}\,\}}}
\newcommand{\lrparent}{\mathop{(\,\mathop{,}\,)}}

\renewcommand{\star}{{}^*}

\newcommand{\updot}{{}^{\bullet}}
\newcommand{\vdual}{{}^{\vee}}

\newcommand\addots{\mathinner{\mkern1mu\raise0pt\vbox{\kern7pt\hbox{.}}\mkern2mu\raise3pt\hbox{.}\mkern2mu\raise6pt\hbox{.}\mkern1mu}}

\newcommand{\restricted}[1]{{}_{\vert}{}_{#1}}
\newcommand{\set}[2]{\{\,#1 \mid #2\,\}}

\newcommand\lto{\longrightarrow}
\newcommand\ltoover[1]{\mathrel{\smash{\overset{#1}{\lto}}}}
\newcommand\varto[1]{\mathrel{\hbox to #1pt{\rightarrowfill}}}
\newcommand\vartoover[2]{\mathrel{{\overset{#2}{\varto{#1}}}}}

\newcommand{\bijective}{\leftrightarrow}

\newcommand{\mono}{\hookrightarrow}

\newcommand{\sends}{\mapsto}

\newcommand{\iso}{\overset{\sim}{\to}}
\newcommand{\liso}{\overset{\sim}{\lto}}

\DeclareMathOperator{\Nilp}{(Nilp)}

\newcommand{\Mss}{\Mcal^{\rm ss}}

\newcounter{segno}[section]
\renewcommand{\thesegno}{\thesection.\arabic{segno}}

\newcounter{thmno}[section]
\renewcommand{\thethmno}{\thesection.\arabic{thmno}}

\newskip\segskipamount
\newskip\procskipamount
\newskip\interskipamount
\newskip\exskipamount
\newskip\refskipamount
\newcommand{\segskip}{\vskip\segskipamount}
\newcommand{\procskip}{\vskip\procskipamount}
\newcommand{\interskip}{\vskip\interskipamount}
\newcommand{\exskip}{\vskip\exskipamount}
\newcommand{\refskip}{\vskip\refskipamount}

\newcommand{\segbreak}{\par
   \ifdim\lastskip<\segskipamount\removelastskip
   \penalty-200
   \segskip\fi}

\newcommand{\procbreak}{\par
   \ifdim\lastskip<\procskipamount\removelastskip
   \penalty-100
   \procskip\fi}

\newcommand{\exbreak}{\par
   \ifdim\lastskip<\exskipamount\removelastskip
   \penalty-100
   \exskip\fi}

\newcommand{\titlebreak}{\par%
\ifdim\lastskip<\interskipamount\removelastskip%
\penalty10000%
\interskip\fi%
\noindent}%

\newcommand{\interbreak}{\par%
\ifdim\lastskip<\interskipamount\removelastskip%
\penalty-100%
\interskip\fi%
\noindent\ignorespaces}%

\newcommand{\introbreak}{\par
   \ifdim\lastskip<\segskipamount\removelastskip
   \nobreak
   \segskip\fi}

\newcommand{\refbreak}{\par
   \ifdim\lastskip<\refskipamount\removelastskip
   \penalty-100
   \refskip\fi}

\numberwithin{equation}{segno}

\newcommand{\marginrule}{\marginpar{\rule[-10.5mm]{1mm}{10mm}}}

\newenvironment{segment}[2]{%
\segbreak%
\refstepcounter{segno}%
\noindent\textbf{(\thesegno)\ }%
\textbf{#1.\ }%
\label{#2}%
\titlebreak\ignorespaces}%
{\par}%

{\par}%

\newenvironment{proclamation}[1]{%
  \procbreak%
  \refstepcounter{thmno}%
  \noindent\textbf{#1\ \thethmno.\ }%
  \em}%
  {\par\interskip}%

\newenvironment{remark}[1]{%
  \procbreak%
  \refstepcounter{thmno}%
  \noindent\textbf{#1\ \thethmno.\ }%
  }%
  {\par\interskip}%

\newenvironment{intermediate}{%
  \segbreak%
  \noindent\ignorespaces}%
  {\par\segskip}%

\newenvironment{proof}{%
\noindent%
{\em Proof.\ }}%
{\hspace*{\fill}$\Box$\\[0.5\topskip]\noindent}

  {\par}%

\newcounter{listcounter}
\newcounter{deflistcounter}

\newcounter{exlistcounter}

\newskip{\itemsepamount}
\newskip{\topsepamount}
\newenvironment{assertionlist}{%
  \begin{list}
    {\upshape (\arabic{listcounter})}
    {\setlength{\leftmargin}{18pt}
     \setlength{\rightmargin}{0pt}
     \setlength{\itemindent}{0pt}
     \setlength{\labelsep}{5pt}
     \setlength{\labelwidth}{13pt}
     \setlength{\listparindent}{\parindent}
     \setlength{\parsep}{0pt}
     \setlength{\itemsep}{\itemsepamount}
     \setlength{\topsep}{\topsepamount}
     \usecounter{listcounter}}}
  {\end{list}}

\newenvironment{definitionlist}{%
  \begin{list}
    {\upshape (\alph{deflistcounter})}
    {\setlength{\leftmargin}{18pt}
     \setlength{\rightmargin}{0pt}
     \setlength{\itemindent}{0pt}
     \setlength{\labelsep}{5pt}
     \setlength{\labelwidth}{13pt}
     \setlength{\listparindent}{\parindent}
     \setlength{\parsep}{0pt}
     \setlength{\itemsep}{\itemsepamount}
     \setlength{\topsep}{\topsepamount}
     \usecounter{deflistcounter}}}
  {\end{list}}

\newenvironment{equivlist}{%
  \begin{list}
    {\upshape (\roman{listcounter})}
    {\setlength{\leftmargin}{18pt}
     \setlength{\rightmargin}{0pt}
     \setlength{\itemindent}{0pt}
     \setlength{\labelsep}{5pt}
     \setlength{\labelwidth}{13pt}
     \setlength{\listparindent}{\parindent}
     \setlength{\parsep}{0pt}
     \setlength{\itemsep}{\itemsepamount}
     \setlength{\topsep}{\topsepamount}
     \usecounter{listcounter}}}
  {\end{list}}

\newenvironment{bulletlist}{%
  \begin{list}
    {\upshape \textbullet}
    {\setlength{\leftmargin}{18pt}
     \setlength{\rightmargin}{0pt}
     \setlength{\itemindent}{0pt}
     \setlength{\labelsep}{6pt}
     \setlength{\labelwidth}{12pt}
     \setlength{\listparindent}{\parindent}
     \setlength{\parsep}{0pt}
     \setlength{\itemsep}{\itemsepamount}
     \setlength{\topsep}{\topsepamount}}}
  {\end{list}}

% ...
\begin{document}

\title{The supersingular locus of the Shimura variety of $\GU(1,n-1)$ II}

\author{Inken Vollaard, Torsten Wedhorn}

\maketitle

%------------------------------------------------------------------

\noindent{\scshape Abstract.\ }
We complete the study of the supersingular locus $\Mss$ in the fiber at $p$ of a Shimura variety attached to a unitary similitude group $\GU(1,n-1)$ over $\QQ$ in the case that $p$ is inert. This was started by the first author in~\cite{Vo_Uni} where complete results were obtained for $n =2,3$. The supersingular locus $\Mss$ is uniformized by a formal scheme $\Ncal$ which is a moduli space of so-called unitary $p$-divisible groups. It depends on the choice of a unitary isocrystal $\Nbf$. We define a stratification of $\Ncal$ indexed by vertices of the Bruhat-Tits building attached to the reductive group of automorphisms of $\Nbf$. We show that the combinatorial behaviour of this stratification is given by the simplicial structure of the building. The closures of the strata (and in particular the irreducible components of $\Ncal_{\rm red}$) are identified with (generalized) Deligne-Lusztig varieties. We show that the Bruhat-Tits stratification is a refinement of the Ekedahl-Oort stratification and also relate the Ekedahl-Oort strata to Deligne-Lusztig varieties. We deduce that $\Mss$ is locally a complete intersection, that its irreducible components and each Ekedahl-Oort stratum in every irreducible component is isomorphic to a Deligne-Lusztig variety, and give formulas for the number of irreducible components of every Ekedahl-Oort stratum of $\Mss$.

%==================================================================

\section*{Introduction}

Let $(G,X)$ be a Shimura datum, where $G$ is a reductive group over $\QQ$ such that $G_{\RR}$ is isomorphic to the group of unitary similitudes of a hermitian space of signature $(1,n-1)$ for some integer $n > 1$. The reflex field of $(G,X)$ is contained in a quadratic imaginary extension $E$ of $\QQ$. We fix a prime $p > 2$ such that $p$ is {\em inert} in $E$ and we assume that $G$ is unramified at $p$. Let $C^p \subset G(\AA^p_f)$ be an open compact subgroup and let $\Mcal_{C^p}$ be the associated moduli space of abelian varieties defined by Kottwitz \cite{Ko_ShFin} for these PEL-Shimura data. It is defined over $\ZZ_{p^2} := \Ocal_{E_v}$ where $v$ is the unique place of $E$ lying over $p$. We refer to Section~\ref{Msupersing} for the precise definition.

In this paper we complete the study of the supersingular locus $\Mcal^{\rm ss}_{C^p}$ of the reduction modulo $p$ of $\Mcal_{C^p}$, which was started in \cite{Vo_Uni}. The uniformization theorem of Rapoport and Zink (\cite{RZ_Period}~Theorem~6.30) provides us with an isomorphism (see~\eqref{NcalMcal}) 
\[
\coprod_j \Gamma_j\backslash \Ncal_{\red} \liso \Mcal_{C^p}^{\rm ss} \tag{U}
\]
where $\Ncal_{\red}$ is the underlying reduced subscheme of a moduli space $\Ncal$ of quasi-isogenies of unitary $p$-divisible groups $(X,\rho_X)$ (see~\eqref{Modp} for its definition) and where $(\Gamma_j)_j$ is a finite family of discrete groups acting continuously on $\Ncal$. The moduli space $\Ncal$ is represented by a separated formal scheme locally formally of finite type over $\Spf(\ZZ_{p^2})$. We study the structure of $\Ncal$ and of $\Ncal_{\red}$. 

By \cite{Vo_Uni}, we have a decomposition $\Ncal = \coprod_{i \in \ZZ} \Ncal_i$ where $\Ncal_i$ classifies \mbox{$p$-iso}genies of height $ni$. Moreover, $\Ncal_i$ is nonempty if and only if $ni$ is even and in this case $\Ncal_i$ is isomorphic to $\Ncal_0$. For the remainder of the introduction let $i$ be fixed such that $ni$ is even.

We define a stratification of $\Ncal_{i,\red}$ by locally closed subschemes $\Ncal^0_{\Lambda}$, where $\Lambda$ runs through the set of vertices of a Bruhat-Tits building $\Bcal$ associated with the data above. More precisely, denote by $\Nbf$ a superspecial unitary isocrystal over $\FF_{p^2}$ (see~\eqref{StandPDiv}). The automorphisms of $\Nbf$ form an algebraic group $J$ over $\QQ_p$ which is an inner form of $G_{\QQ_p}$. Let $\Jtilde$ be the derived group of $J$, denote by $\Bcal = \Bcal(\Jtilde,\QQ_p)$ its simplicial Bruhat Tits building, and let $\Bcal_0$ be the set of vertices of $\Bcal$. Every such vertex can be considered as a lattice in a nondegenerate $(K/\QQ_p))$-hermitian space, where $K$ is a quadratic unramified extension of $\QQ_p$. For each such lattice $\Lambda \in \Bcal_0$ we define a projective subscheme $\Ncal_{\Lambda}$ of $\Ncal_{i,\red}$ (see \eqref{VertexScheme} for the definition).

In \cite{Vo_Uni} subsets $\Vcal(\Lambda)(k)$ of $\Ncal_i(k)$ were defined for each algebraically closed extension $k$ of $\FF_{p^2}$ and it was conjectured that these subsets are the sets of a $k$-valued points of a subscheme of $\Ncal_i$. This conjecture was proved in loc.~cit.\ for $n = 2$ and $n = 3$. As $\Ncal_{\Lambda}(k) = \Vcal(\Lambda)(k)$, the construction of $\Ncal_{\Lambda}$ proves this conjecture for arbitrary $n$.

To describe the structure of $\Ncal_{\Lambda}$ we introduce the following notation. For each $\Lambda \in \Bcal_0$ let $\Jtilde_{\Lambda}$ be the maximal reductive quotient of the special fibre of the parahoric group scheme of $\Jtilde$ associated with the vertex $\Lambda$. Then $\Jtilde_{\Lambda}$ is the special unitary group of an $(\FF_{p^2}/\FF_p)$-hermitian space $V_{\Lambda}$ of odd dimension, say $2t(\Lambda) + 1$. The group $\Jtilde(\QQ_p)$ acts on $\Bcal_0$ and the map $\Lambda \sends t(\Lambda)$ induces a bijection between the set of $\Jtilde(\QQ_p)$-orbits on $\Bcal_0$ and the set of integers $\{0,1,\dots,[(n-1)/2]\}$. We call $t(\Lambda)$ the orbit type of $\Lambda \in \Bcal_0$.

The schemes $\Ncal_{\Lambda}$ have the following properties (see Corollary~\ref{PropNLambda}).

\procskip

\noindent{\bf Theorem A.}\ {\em $\Ncal_{\Lambda}$ is projective, smooth, and geometrically irreducible of dimension $t(\Lambda)$. If $t(\Lambda) = t(\Lambda')$, the $\FF_{p^2}$-schemes $\Ncal_{\Lambda}$ and $\Ncal_{\Lambda'}$ are isomorphic.
}

\procskip

For the proof of Theorem~A we construct an isomorphism of $\Ncal_{\Lambda}$ with a Deligne-Lusztig variety $Y_{\Lambda}$ of the group $\Jtilde_{\Lambda}$ (see~\eqref{YLambda}). Here we mean by a Deligne-Lusztig variety of an algebraic group $H$ (defined over a finite field) a variety classifying parabolic subgroups $P$ of $H$ of a fixed conjugacy type such that $P$ and $F(P)$ are in a fixed relative position, where $F$ is the Frobenius. This generalizes the varieties considered by Deligne and Lusztig who studied the case where the conjugacy class of parabolic subgroups is the class of Borel subgroups. In~\eqref{DLVar} we collect some results about these (generalized) Deligne-Lusztig varieties.

The intersection behaviour of the subschemes $\Ncal_{\Lambda}$ is given by the following result (see Theorem~\ref{CombNcalLambda} and Proposition~\ref{Ncal0Open}).

\procskip

\noindent{\bf Theorem B.}\ {\em Let $\Lambda, \Lambda' \in \Bcal_0$.
\begin{assertionlist}
\item
The Bruhat-Tits stratum $\Ncal_{\Lambda'}$ is contained in $\Ncal_{\Lambda}$ if and only if $\Lambda' \subset \Lambda$.
\item
Two Bruhat-Tits strata $\Ncal_{\Lambda}$ and $\Ncal_{\Lambda'}$ have a nonempty interscetion if and only if $\Lambda \cap \Lambda' \in \Bcal_0$, and in this case $\Ncal_{\Lambda} \cap \Ncal_{\Lambda'} = \Ncal_{\Lambda \cap \Lambda'}$ scheme-theoretically.
\item
Let $\Ncal_{\Lambda}^0$ be the open subscheme of $\Ncal_{\Lambda}$ which is the complement of the union of the (finitely many) closed subschemes $\Ncal_{\Lambda'}$ with $\Ncal_{\Lambda'} \subsetneq \Ncal_{\Lambda}$. Then the natural morphism
\[
\coprod_{\Lambda \in \Bcal_0}  \Ncal_{\Lambda}^0 \lto \Ncal_i
\]
is bijective.
\item
The subscheme $\Ncal_{\Lambda}^0$ is open and dense in $\Ncal_{\Lambda}$.
\end{assertionlist}
}

\procskip

In particular, we obtain the stratification of $\Ncal_i$ by the $\Ncal^0_{\Lambda}$ mentioned above. Theorem~B shows that the intersection behaviour and the inclusion relation of the subschemes $\Ncal_{\Lambda}$ are given by the combinatorial structure of the Bruhat-Tits building of $\Jtilde$. Therefore we call this stratification the \emph{Bruhat-Tits stratification}. As an application we show that any kind of intersection behaviour within the bounds given by Theorem~B can occur (Proposition~\ref{InterArb}) and we give formulas how many closed Bruhat-Tits strata of a fixed dimension are contained in and are containing a given one (Corollary~\ref{NLLPrime}).

Moreover, we obtain the following corollary (Theorem~\ref{NLambdaConn} and Theorem~\ref{LocalNcal}).

\procskip

\noindent{\bf Corollary C.}\ {\em $\Ncal_i$ is geometrically connected of pure dimension $[(n-1)/2]$ and locally a complete intersection. The irreducible components of $\Ncal_{i,\red}$ are precisely the projective subschemes $\Ncal_{\Lambda}$ for $\Lambda \in \Bcal_0$ with $t(\Lambda) = [(n-1)/2]$.}

\procskip

Then we compare the Bruhat-Tits stratification with the Ekedahl-Oort stratification of $\Ncal_i$ given by the isomorphism class of the $p$-torsion of the unitary $p$-divisible group (defined in~\eqref{EOStrata}). We show that the Ekedahl-Oort strata are smooth (Theorem~\ref{EOSmooth}) and that the former stratification is a refinement of the latter one (see~\eqref{EONLambda}):

\procskip

\noindent{\bf Theorem D.}\ {\em Let $x,x' \in \Ncal_0(k)$, where $k$ is an algebraically closed extension of $\FF_{p^2}$. Then these points lie in the same Ekedahl-Oort stratum if and only $x \in \Ncal_{\Lambda}^0(k)$ and $x' \in \Ncal_{\Lambda'}^0(k)$ with $t(\Lambda) = t(\Lambda')$.}

\procskip

We obtain the following corollary (see Corollary~\ref{EODL} and Theorem~\ref{LocalNcal}).

\procskip

\noindent{\bf Corollary E.}\ {\em The smooth locus of $\Ncal_i$ is the open Ekedahl-Oort stratum of $\Ncal_i$. Each Ekedahl-Oort stratum of $\Ncal_{\Lambda}$ is isomorphic to a Deligne-Lusztig variety.}

\procskip

Using the theory of $p$-adic uniformization mentioned above we obtain now the following result for the local structure of the supersingular locus $\Mss_{C^p}$ (see Theorem~\ref{NMss}).

\procskip

\noindent{\bf Theorem F.}\ {\em The supersingular locus $\Mss_{C^p}$ is of pure dimension $[(n-1)/2]$ and locally of complete intersection. Its smooth locus is the open Ekedahl-Oort stratum $\Mss_{C^p}([(n-1)/2])$.}

\procskip

Note that by~\cite{Wd_OS} all Ekedahl-Oort strata of $\Mcal_{C^p}$ are smooth. We also obtain formulas for the number of connected components of $\Mss_{C^p}$ and the number of irreducible components of all Ekedahl-Oort strata (and in particular for the number of irreducible components for $\Mss_{C^p}$) (see Proposition~\ref{IrredMss} and Proposition~\ref{ConnMss}). These formulas in particular show that within a connected component of $\Mss_{C^p}$ the number of irreducible components of a given Ekedahl-Oort stratum becomes arbitrarily large when $p$ goes to infinity.

If $C^p$ is sufficiently small, the irreducible components of $\Mss_{C^p}$ are isomorphic to the irreducible components of $\Ncal_i$ (which are all pairwise isomorphic by Theorem~A). We illustrate the above results in the following low-dimensional example (the case $n=3$ has already been studied in~\cite{Vo_Uni}, for $n = 4$ see~\eqref{ExGlobal} and Proposition~\ref{BTIsom}).

\procskip

\noindent{\bf Example G.}\ {\em Let $n=3$ or $n=4$ and let $C^p$ be sufficiently small. The supersingular locus $\Mss_{C^p}$ is projective, equi-dimensional of dimension $1$, and locally a complete intersection. Each irreducible component is isomorphic to the Fermat curve in $\PP^2$ given by the equation $x_0^{p+1} + x_1^{p+1} + x_2^{p+1} = 0$. Its singular points are the superspecial points of $\Mss_{C^p}$ and there are $p^3+1$ of them on any irreducible component. Each superspecial point is the pairwise transversal intersection of $p+1$ (for $n=3$), resp.~of $p^3+1$ (for $n = 4$) irreducible components.} 

\bigskip

The geometric results obtained in this paper are used by S.~Kudla and M.~Rapoport in \cite{KR_CycUniI} where they study the local intersection theory of special cycles in the unramified case.

S.~Harashita has obtained in \cite{Ha_EODL} a description in terms of Deligne-Lusztig varieties for the union of certain Ekedahl-Oort strata for the moduli space of principally polarized abelian varieties. In fact, some results of loc.~cit.\ can be reformulated using the language of Bruhat-Tits buildings and this reformulation initiated in part this paper. Harashita's results have been refined by M.~Hoeve in \cite{Ho_EOSuperSing} using so-called fine Deligne-Lusztig varieties. We remark that this is analogous to the results described here, as the Deligne-Lusztig varieties of Corollary~\ref{EODL} can be also described as fine Deligne-Lusztig varieties.

\bigskip

The paper is organized as follows. The first section fixes notations and defines the formal scheme $\Ncal$. Section~2 studies the Ekedahl Oort stratification of $\Ncal$. In Section~3 the schemes $\Ncal_{\Lambda}$ are defined and it is proved that they are isomorphic to certain Deligne-Lusztig varieties. This is used in Section~4 to prove the main results on the structure of $\Ncal$. In Section~5 we explain how these results are applied to the structure of the supersingular locus of Shimura varieties attached to certain unitary groups.

\bigskip

\noindent\textsc{Acknowledgments}: We are very grateful to M.~Rapoport whose questions motivated this paper and whose comments were very helpful, and to E.~Lau for a helpful discussion. The second author also gratefully acknowledges the hospitality of the Erwin-Schrödinger Institute in Vienna in 2007, where this work was started.

%------------------------------------------------------------------

\section{The Moduli space $\Nscr$ of $p$-divisible groups}

\begin{segment}{The local PEL-datum}{PEL}
Let $p > 2$ be a fixed prime number. Let $K$ be an unramified extension of $\QQ_p$ of degree $2$ and denote by $\star$ the nontrivial Galois automorphism of $K$ over $\QQ_p$. Let $V$ be a finite-dimensional $K$-vector space and set $n = \dim_K(V)$. Let $\lrangle$ be a $\QQ_p$-valued skew-hermitian form on $V$, i.e., $\lrangle\colon V \times V \to \QQ_p$ is an alternating $\QQ_p$-bilinear form such that $\langle av, w\rangle = \langle v,a\star w\rangle$. Let $G$ be the algebraic group over $\QQ_p$ such that
\[
G(R) = \set{g \in \GL_{K \otimes R}(V_R)}{\exists c \in R^{\times}\colon \langle gv, gw \rangle = c \langle v,w \rangle \forall v,w \in V \otimes R}
\]
for every $\QQ_p$-algebra $R$. This is a reductive group over $\QQ_p$.

There exists a unique skew-hermitian form $(\ ,\ )'\colon V \times V \to K$ such that $\langle v,w \rangle = \Tr_{K/\QQ_p}(v,w)'$. Moreover, if $\delta \in K\cross$ with $\delta\star = -\delta$, $(v,w) := \delta(v,w)'$ is a hermitian form, and $G$ is the reductive group of similitudes of the hermitian space $(V, (\ ,\ ))$.

Let $O_K$ be the ring of integers of $K$. We assume that there exists an \mbox{$O_K$-lattice} $\Gamma$ such that $\lrangle$ induces a perfect $\ZZ_p$-pairing on $\Gamma$. This implies that $G$ has a reductive model over $\ZZ_p$, namely the group of $O_K$-linear symplectic similitudes of $(\Gamma, \lrangle\restricted{\Gamma \times \Gamma})$.
\end{segment}

\begin{segment}{Unitary $p$-divisible groups}{UnitaryPDiv}
Let $\FF_{p^2}$ be a finite field with $p^2$ elements and define $\ZZ_{p^2} := W(\FF_{p^2})$ and let $\QQ_{p^2}$ be its field of fractions. We denote by $\varphi_0$ and $\varphi_1 = \varphi_0 \circ \star$ the two $\QQ_p$-isomorphisms of $K$ to $\QQ_{p^2}$. Let $\Nilp$ be the category of $\ZZ_{p^2}$-schemes $S$, such that $p$ is locally nilpotent on $S$.

We fix an integer $r$ with $0 \leq r \leq n$. For each $S \in \Nilp$ a {\em $p$-divisible $O_K$-module of signature $(r,n-r)$ over $S$} is a pair $X = (X,\iota_X)$, where $X$ is a $p$-divisible group over $S$ with an $O_K$-action $\iota_X\colon O_K \to \End(X)$ satisfying
\begin{equation}\label{signature}
{\rm charpol}(\iota(a)|\Lie(X)) = (T - \varphi_0(a))^r(T - \varphi_1(a))^{n-r} \in \ZZ_{p^2}[T]
\end{equation}
for all $a \in \Ocal_K$. We call this the {\em signature condition $(r,n-r)$}.
Note that \eqref{signature} implies $\rk_{\Oscr_S}(\Lie(X)) = n$.

We denote by ${}^tX$ the dual $p$-divisible group endowed with the $O_K$-action $\iota_{{}^tX}(a) = {}^t(\iota_X(a^{*}))$. A {\em unitary $p$-divisible group of signature $(r,n-r)$} is a triple $X = (X,\iota_X,\lambda_X)$ where $(X,\iota_X)$ is a $p$-divisible $O_K$-module of signature $(r,n-r)$ and where $\lambda_X\colon X \to {}^tX$ is a $p$-principal $O_K$-linear polarization. The existence of $\lambda_X$ implies
\[
\height(X) = \rk_{\Oscr_S}(\Lie(X)) + \rk_{\Oscr_S}(\Lie({}^tX)) = 2n.
\]

We call two unitary $p$-divisible groups $X$ and $Y$ {\em isomorphic} if there exists an $\Ocal_K$-linear isomorphism $\alpha\colon X \to Y$ such that ${}^t\alpha \circ \lambda_Y \circ \alpha$ is a $\ZZ\cross_p$-multiple of $\lambda_X$.

Finally, by repeating the definitions above for at level 1 truncated Barsotti-Tate groups we obtain the notion of a {\em unitary $\BT_1$ of signature $(r,n-r)$}.
\end{segment}

\begin{segment}{Dieudonn\'e modules of unitary $p$-divisible groups}{DieudUni}
Now assume that $S = \Spec(k)$ where $k$ is a perfect field extension of $\FF_{p^2}$. Then the covariant Dieudonn\'e module of a unitary $p$-divisible groups of signature $(r,n-r)$ is a Dieudonn\'e module $(M,\Fcal,\Vcal)$ of height $2n$ over $k$ endowed with an $\Ocal_K$-action and a perfect alternating $W(k)$-bilinear pairing $\lrangle$ on $M$ satisfying
\begin{equation}\label{OKcomp}
\langle \Fcal x, y \rangle = \langle x, \Vcal y\rangle^{\sigma}, \qquad \langle ax, y \rangle = \langle x, a\star y \rangle
\end{equation}
for all $x,y \in M$ and $a \in \Ocal_K$. Here $\sigma$ denotes the Frobenius on $W(k)$.

Via the decomposition
\begin{equation}\label{Decomp01}
\Ocal_K \otimes_{\ZZ_p} W(k) \liso W(k) \times W(k), \quad a \otimes w \sends (\varphi_0(a)w, \varphi_1(a)w),
\end{equation}
the $\Ocal_K$-action on $M$ defines a decomposition $M = M_0 \oplus M_1$ and the conditions \eqref{OKcomp} mean that $\Fcal$ and $\Vcal$ are homogeneous of degree $1$ with respect to this decomposition and that $M_0$ and $M_1$ are totally isotropic with respect to $\lrangle$. Finally the signature condition~\eqref{signature} is equivalent to the equalities
\[
\dim_k(M_0/\Vcal M_1) = r, \qquad \dim_k(M_1/\Vcal M_0) = n - r.
\]
We call such tuples $M = (M,\Fcal,\Vcal, M = M_0 \oplus M_1, \lrangle)$ {\em unitary Dieudonn\'e modules of signature $(r,n-r)$}.

If $M = (M,\Fcal,\Vcal, M = M_0 \oplus M_1, \lrangle)$ is a unitary Dieudonn\'e module of signature $(r,n-r)$ over $k$ we can consider its reduction $\Mbar = M \otimes_{W(k)} k$ which is a $\ZZ/2\ZZ$-graded $k$-vector space of dimension $2n$ endowed with $\Fcal$ and $\Vcal$ of degree $1$ and a perfect alternating pairing $\lrangle$ such that $\Mbar_0$ and $\Mbar_1$ are totally isotropic with respect to $\lrangle$ and such that \mbox{$\dim_k(\Mbar_0/\Vcal\Mbar_1) = r$} and \mbox{$\dim_k(\Mbar_1/\Vcal\Mbar_0) = n-r$}. We call such a tuple
\[
\Mbar = (\Mbar,\Fcal,\Vcal, \Mbar = \Mbar_0 \oplus \Mbar_1, \lrangle)
\]
a {\em unitary Dieudonn\'e space of signature $(r,n-r)$}. 
\end{segment}

\begin{segment}{The standard unitary $p$-divisible group $\Xbf$}{StandPDiv}
We fix a unitary $p$-divisible group $\Xbf = (\Xbf,\igbf,\lgbf)$ over $\FF_{p^2}$ of signature \mbox{$(r,n-r)$}. We let $(\Mbf,\Fbf,\Vbf)$ be the covariant Dieudonn\'e module of $\Xbf$. Then $\Mbf$ is a free $\ZZ_{p^2}$-module of rank $2n$. Via $\igbf$ it is endowed with an $O_K$-action and $\lgbf$ defines a perfect symplectic pairing on $\Mbf$. Therefore $\Mbf$ is a unitary Dieudonn\'e module over $\FF_{p^2}$.

We denote by $(\Nbf,\Fbf) = (\Mbf,\Fbf) \otimes \QQ_{p^2}$ the isocrystal of $\Xbf$. 

We assume that $\Xbf$ is superspecial and that the isocrystal $(\Nbf,\Fbf)$ of $\Xbf$ is generated by elements $n$ such that $\Fbf^2n = pn$. In particular, $\Nbf$ is decent in the sense of \cite{RZ_Period}~2.13. As $\Fbf^2$ is a $\QQ_{p^2}$-linear map, we have $\Fbf^2 = p\id_{\Nbf}$ and therefore $\Fbf = \Vbf$.

Such a triple $(\Xbf,\igbf,\lgbf)$ exists. Indeed, we define its unitary Dieudonn\'e module $\Mbf$ as follows. Fix $\delta \in \ZZ_{p^2}^{\times}$ such that $\sigma(\delta) = -\delta$. We consider $\SStilde = O_K \otimes \ZZ_{p^2} = \ZZ_{p^2} \times \ZZ_{p^2}$ and set $g := (1,0), h := (0,1) \in \SStilde$. Define a $(\ZZ_{p^2},\sigma)$-linear map $\Fbf$ on $\SStilde$ by $\Fbf(h) = g$ and $\Fbf(g) = ph$ and define a $(\ZZ_{p^2},\sigma^{-1})$-linear map $\Vbf$ on $\SStilde$ by $\Vbf(h) = g$ and $\Vbf(g) = ph$. We define a perfect $\ZZ_{p^2}$-pairing on $\SStilde$ by $\langle g,h\rangle = \delta$. This makes $\SStilde$ into a unitary Dieudonn\'e module of signature $(0,1)$ and we set
\[
\Mbf = \SStilde^{n-r} \oplus \sigma^*(\SStilde)^r,
\]
where $\sigma^*(\SStilde)$ is $\SStilde$ as an $O_K \otimes \ZZ_{p^2}$-module with symplectic pairing but where for the definition of $\Fbf$ and $\Vbf$ the roles of $g$ and $h$ are interchanged.
\end{segment}

\begin{segment}{The moduli space of $p$-divisible groups}{Modp}
For each scheme $S \in \Nilp$ we set $\Sbar = S \otimes_{\ZZ_{p^2}} \FF_{p^2}$. We define a set-valued functor $\Ncal$ on $\Nilp$ which sends a scheme $S \in \Nilp$ to the set of isomorphism classes of tuples $X = (X,\iota_X,\lambda_X,\rho_X)$. Here $(X,\iota_X,\lambda_X)$ is a unitary \mbox{$p$-divisible} group over $S$ of signature $(r,n-r)$ and $\rho_X$ is an $\Ocal_K$-linear quasi-isogeny
\[
\rho_X\colon X \times_S \Sbar  \lto \Xbf \times_{\Spec(\FF_{p^2})} \Sbar
\]
such that ${}^t\rho_X \circ \lambda_X \circ \rho_X$ is a $\QQ\cross_p$-multiple of $\lambda_X$ in $\Hom_{\Ocal_K}(X,{}^tX) \otimes_{\ZZ_p} \QQ_p$.

By~\cite{RZ_Period}~Corollary~3.40 the functor $\Ncal$ is representable by a formal scheme over $\Spf(\ZZ_{p^2})$ which is locally formally of finite type. It follows from~loc.~cit.\ Proposition~2.32 that every irreducible component of $\Ncal_{\red}$ is projective, in particular $\Ncal$ is separated.

If $X$ is an $S$-valued point of $\Ncal$, the height of $\rho_X$ (considered as locally constant function on $S$) is a multiple of $n$ by \cite{Vo_Uni}~1.7 and we obtain a decomposition
\[
\Ncal = \coprod_{i \in \ZZ} \Ncal_i,
\]
where $\Ncal_i$ is the open and closed formal subscheme of $\Ncal$ where $\rho_X$ has height $ni$. Then by \cite{Vo_Uni}~1.8 and~1.22 we have
\begin{equation}\label{NcalEmpty}
\Ncal_i \ne \emptyset \iff \text{$ni$ is even}.
\end{equation}
Note that in this case we have $\Ncal_i \cong \Ncal_0$ by Proposition~\ref{JTransConn} below.
\end{segment}

\begin{segment}{The group $J$ of automorphisms of $\Xbf$}{GroupJ}
Let $J$ be the algebraic group of automorphisms of the unitary isocrystal $\Nbf$, i.e., for any $\QQ_p$-algebra $R$ we denote by $J(R) = J_{\Nbf}(R)$ the group of $(K \otimes_{\QQ_p} \QQ_{p^2}) \otimes_{\QQ_p} R$-linear symplectic similitudes of $\Nbf \otimes R$ which commute with $\Fbf \otimes \id_R$. This is a reductive group over $\QQ_p$ which is an inner form of $G$. By Dieudonn\'e theory, $J(\QQ_p)$ is nothing but the group of quasi-isogenies $\Xbf \to \Xbf$ and therefore $J(\QQ_p)$ acts on $\Ncal$.

As we assumed that $\Nbf$ is generated by elements $n$ such that $\Fbf^2n = pn$, the arguments of \cite{Vo_Uni}~1.19 show that if $k$ is any perfect field extension of $\FF_{p^2}$ the group of automorphisms of the unitary isocrystal associated with $\Xbf \otimes_{\FF_{p^2}} k$ is again $J$.

We fix an element $\delta \in \ZZ_{p^2}^{\times}$ such that $\delta^{\sigma} = -\delta$ and define a nondegenerate $\sigma$-hermitian form on the $\QQ_{p^2}$-vector space $\Nbf_0$ by
\begin{equation}\label{Curly}
\{x,y\} := \delta\langle x, Fy \rangle.
\end{equation}
The $K$-linearity of $g \in J(R)$ implies that $g(\Nbf_0 \otimes R) = \Nbf_0 \otimes R$ and restricting $g$ to $\Nbf_0 \otimes R$ defines an isomorphism of $J$ with the algebraic group $\GU(\Nbf_0,\lrcurly)$ of unitary similitudes of the hermitian space $(\Nbf_0,\lrcurly)$.

We set
\begin{equation}\label{MatrixHerm}
T_{\rm odd} := \begin{pmatrix}
& & & 1 \\
& & 1 \\
& \addots \\
1
\end{pmatrix}, \qquad
T_{\rm even} := \begin{pmatrix}
& & & p \\
& & 1 \\
& \addots \\
1
\end{pmatrix}.
\end{equation}
By \cite{Vo_Uni}~Proposition~1.18 (see also the proof of loc.~cit.\ Lemma~1.20) there exists a $\QQ_{p^2}$-basis of $\Nbf_0$ such that the hermitian form $\lrcurly$ is given by $T_{\rm odd}$ if $n$ is odd and by $T_{\rm even}$ if $n$ is even.

Finally, we recall from \cite{Vo_Uni}~1.22:

\begin{proclamation}{Proposition}\label{JTransConn}
Let $i \in \ZZ$ be an integer such that $ni$ is even (i.e.~$\Ncal_i \ne \emptyset$). Then there exists a $g \in J(\QQ_p)$ such that $g(\Ncal_i) = \Ncal_0$. In particular, $\Ncal_i$ is isomorphic to $\Ncal_0$.
\end{proclamation}
\end{segment}

\begin{segment}{Description of the points of $\Ncal$}{Points}
Let $k$ be a perfect field extension of $\FF_{p^2}$. Denote by $\Nbf_k$ the unitary isocrystal $\Nbf \otimes_{\ZZ_{p^2}} W(k)$. Then the remarks in~\eqref{DieudUni} show that by covariant Dieudonn\'e theory we obtain a bijection between the set $\Ncal_i(k)$ and the set of \mbox{$W(k)$-lat}\-tices $M$ in $\Nbf_k$ such that
\begin{bulletlist}
\item
$M$ is stable under $\Fcal$ and $\Vcal$,
\item
$M = M_0 \oplus M_1$ where $M_j = M \cap (\Nbf_j)_k$,
\item
$\dim_k(M_0/\Vcal M_1) = r$ and $\dim_k(M_1/\Vcal M_0) = n - r$,
\item
$M = p^iM^{\perp}$, where $M^{\perp} = \set{x \in \Nbf_k}{\langle x, M \rangle \subset W(k)}$.
\end{bulletlist}
\end{segment}

%------------------------------------------------------------------

\section{EO-stratification of $\Ncal$}

%\begin{intermediate}
From now on we assume that $r = 1$.
%\end{intermediate}

\begin{segment}{Examples of unitary Dieudonn\'e spaces}{Braid}
We give two examples $\SS$ and $\BB(d)$ ($d \geq 1$ an integer) of unitary Dieudonn\'e spaces over $\FF_{p^2}$ which will play an important role in the Ekedahl-Oort stratification of $\Ncal$.

The underlying $\FF_{p^2}$-vector space of $\SS$ is the space of dimension $2$ generated by two elements $g$ and $h$. We set $M_0 = \FF_{p^2}g$ and $M_1 = \FF_{p^2}h$. The alternating form is given by $\langle h,g \rangle = 1$, and we have $\Fcal(h) = -g$ and $\Vcal(g) = h$. Then $\SS$ is a unitary Dieudonn\'e space of signature $(0,1)$. We indicate the definition of $\SS$ by the following diagram
\[\xymatrix{
g \\ h \ar@{-->}@/^/[u]^{-} \ar@/_/[u],
}\]
where $\Fcal$ is given by the broken arrow and $\Vcal$ by the solid arrow. The vectors generating $M_0$ are given in the first line and the vectors generating $M_1$ are given in the second line. Moreover $\lrangle$ induces a perfect pairing on the space generated by the two vectors in the same column.

For $\BB(d)$ the underlying graded $\FF_{p^2}$-vector space $M = M_0 \oplus M_1$ is given by
\[
M_0 = \FF_{p^2}e_1 \oplus \ldots \oplus \FF_{p^2}e_d, \qquad M_1 = \FF_{p^2}f_1 \oplus \ldots \oplus \FF_{p^2}f_d.
\]
The alternating form is defined by
\[
\langle e_i,f_j \rangle = (-1)^{i}\delta_{ij}.
\]
Finally $\Fcal$ and $\Vcal$ are given by
\begin{align*}
\Vcal(f_i) &= e_{i+1},\qquad\text{for $i=1,\ldots,d-1$,} \\
\Vcal(e_d) &= f_1, \\
\Fcal(f_i) &= e_{i-1},\qquad\text{for $i = 2,\ldots,d$,} \\
\Fcal(e_1) &= (-1)^df_d.
\end{align*}
This is a unitary Dieudonn\'e module of signature $(1,d-1)$. With the same convention as above $\BB(d)$ is given by the diagram
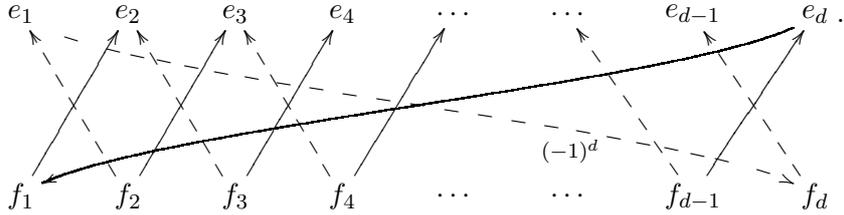
\[\xymatrix{
{e_1} \ar@{-->}@(dr,ul)[ddrrrrrrr]_(.65){(-1)^d} & {e_2} & {e_3} & {e_4} & \cdots & \cdots & {e_{d-1}} & {e_d} \ar@(dl,ur)[ddlllllll] \\
& & & & & & & \\
{f_1} \ar[uur] & {f_2} \ar[uur] \ar@{-->}[uul] & {f_3} \ar[uur] \ar@{-->}[uul] & {f_4} \ar[uur] \ar@{-->}[uul] & \cdots & \cdots & {f_{d-1}} \ar[uur] \ar@{-->}[uul] & f_d \ar@{-->}[uul]
}.\]
With these definitions $\SS$ and $\BB(d)$ are unitary Dieudonn\'e spaces of signature $(0,1)$ and $(1,d-1)$, respectively.
\end{segment}

\begin{segment}{Classification of unitary Dieudonn\'e spaces}{ClassDieud}
We call a unitary Dieudonn\'e module $M$ over a perfect field $k$ {\em supersingular}, if its isocrystal $(M[\frac{1}{p}],\Fcal)$ is supersingular.

By \cite{BW_U}~(4.2) and~(3.6) we have the following classification result.

\begin{proclamation}{Theorem}\label{ClassSupersing}
Let $k$ be algebraically closed and let $\Mbar$ be a unitary Dieudonn\'e space of signature $(1,n-1)$. The following assertions are equivalent.
\begin{assertionlist}
\item
There exists a supersingular unitary Dieudonn\'e module $M$ such that $\Mbar \cong M \otimes_{W(k)} k$.
\item
Every unitary Dieudonn\'e module $M$ such that $\Mbar \cong M \otimes_{W(k)} k$ is supersingular.
\item
There exists a (necessarily unique) integer $0 \leq \sigma \leq (n-1)/2$ such that
\[
\Mbar \cong (\BB(2\sigma+1) \oplus \SS^{n-(2\sigma+1)}) \otimes_{\FF_{p^2}} k,
\]
where $\BB(2\sigma+1)$ and $\SS$ are the unitary Dieudonn\'e spaces defined in~\eqref{Braid}.
\end{assertionlist}
\end{proclamation}

We call $\Mbar$ {\em supersingular} if it satisfies the equivalent conditions of Theorem~\ref{ClassSupersing}. We set
\begin{equation}\label{MMSigma}
\MM_{\sigma} := \BB(2\sigma+1) \oplus \SS^{n-(2\sigma+1)}
\end{equation}
and we denote by $\XXbar_{\sigma}$ the corresponding at level 1 truncated unitary Barsotti-Tate group.
\end{segment}

\begin{segment}{Ekedahl-Oort strata of $\Ncal \otimes_{\ZZ_{p^2}} \FF_{p^2}$}{EOStrata}
Fix an integer $0 \leq \sigma \leq (n-1)/2$. The Ekedahl-Oort stratum $\Ncal(\sigma)$ in $\Ncal \otimes_{\ZZ_{p^2}} \FF_{p^2}$ is defined as follows. For any $\FF_{p^2}$-scheme $S$ we define $\Ncal(\sigma)(S)$ to be the $X \in \Ncal(S)$ such that $X[p]$ is locally for the fppf-topology isomorphic to $(\XXbar_{\sigma})_S$. 

By \cite{Wd_OS}~(6.7) the monomorphism $\Ncal(\sigma) \mono \Ncal \otimes_{\ZZ_{p^2}} \FF_{p^2}$ is representable by an immersion. In particular $\Ncal(\sigma)$ is a formal scheme over $\FF_{p^2}$, locally formally of finite type. These formal schemes are called the {\em Ekedahl-Oort strata} of $\Ncal \otimes_{\ZZ_{p^2}} \FF_{p^2}$. For any perfect extension $k$ of $\FF_{p^2}$ we have
\[
\Ncal(k) = \biguplus_{0 \leq \sigma \leq \frac{n-1}{2}} \Ncal(\sigma)(k).
\]
by Theorem~\ref{ClassSupersing}. 

\begin{proclamation}{Theorem}\label{EOSmooth}
The Ekedahl-Oort strata $\Ncal(\sigma)$ are formally smooth over $\FF_{p^2}$.
% and the immersion $\Ncal(\sigma) \mono \Ncal \otimes_{\ZZ_{p^2}} \FF_{p^2}$ is regular of codimension $[(n-1)/2] - \sigma$.
\end{proclamation}

\begin{proof}
Let $R$ be an $\FF_{p^2}$-algebra. Associating to $R$-valued points $(X,\rho_X)$ of $\Ncal$ its unitary $\BT_1$ $X[p]$ defines a morphism
\[
\alpha\colon \Ncal \otimes \FF_{p^2} \to \Bscr,
\]
where $\Bscr$ is the algebraic stack that classifies unitary $\BT_1$ of signature \mbox{$(1,n-1)$} (cf.~\cite{Wd_OS}~(1.8)). In $\Bscr$ there is the smooth locally closed substack $\Bscr(\sigma)$ classifying unitary $\BT_1$ which are locally for the fppf-topology isomorphic to $\XXbar_{\sigma}[p]$. By definition $\Ncal(\sigma)$ is the inverse image of $\Bscr(\sigma)$ under $\alpha$. Therefore it suffices to show that $\alpha$ is smooth.

By~\cite{Wd_OS}~(2.17), associating to a unitary $p$-divisible group $X$ its $p$-torsion $X[p]$ is formally smooth. Therefore the smoothness of $\alpha$ follows from Drinfeld's theorem that quasi-isogenies can be always deformed uniquely.
\end{proof}
\end{segment}

\begin{segment}{Superspecial gap}{ssgap}
We will now connect the Ekedahl-Oort stratification with the notion of type defined in \cite{Vo_Uni}. 

We fix an integer $i \in \ZZ$ such that $ni$ is even~\eqref{NcalEmpty}. Let $k$ be a perfect field extension of $\FF_{p^2}$ and let $x \in \Ncal_i(k)$. Let $M$ be the unitary Dieudonn\'e module over $k$ of signature $(1,n-1)$ corresponding to $x$~\eqref{Points} and let $N = M \otimes_{\ZZ} \QQ$ be its isocrystal. We set $\tau := p^{-1}\Fcal^2\colon N \to N$ and define
\[
\Lambda^+(M) := \sum_{l \geq 0} \tau^l(M), \qquad \Lambda^-(M) := \bigcap_{l \geq 0} \tau^l(M).
\]
These are $\Ocal_K \otimes W(k)$-submodules and therefore there is a decomposition \mbox{$\Lambda^+(M) = \Lambda^+(M)_0 \oplus \Lambda^+(M)_1$}. Note that $\Lambda^+(M)_j$ is $\tau$-invariant for $j = 0,1$. Similar statements hold for $\Lambda^-(M)$. We set
\begin{equation}\label{DefSigma}
\begin{aligned}
\sigma^+(M_j) &:= \inf\set{d \geq 0}{\Lambda^+(M)_j = \sum_{l = 0}^d \tau^l(M_j)},\\
\sigma^-(M_j) &:= \inf\set{d \geq 0}{\Lambda^-(M)_j = \bigcap_{l = 0}^d \tau^l(M_j)}.
\end{aligned}
\end{equation}

The pairing $p^{-i}\lrangle$ induces a perfect pairing on $M$ and therefore a perfect duality between $\Lambda^-(M)_0$ and $\Lambda^+(M)_1$ and between $\Lambda^-(M)_1$ and $\Lambda^+(M)_0$. Hence
\begin{equation}\label{SigmaProp}
\begin{aligned}
\sigma^+(M_0) = \sigma^-(M_1) \leq \frac{n-1}{2}, &\qquad \sigma^+(M_1) = \sigma^-(M_0) \leq \frac{n+1}{2},\\
\sigma^+(M_j) = [\Lambda^+(M_j) : M_j], &\qquad \sigma^-(M_j) = [M_j : \Lambda^-(M_j)],\\
p\Lambda^+(M) \subset \Lambda^-(M) \subset &M \subset \Lambda^+(M) \subset p^{-1}\Lambda^-(M).
\end{aligned}
\end{equation}
by \cite{Vo_Uni}~Lemma~2.2. Moreover we have
\[
\sigma^+(M_0) = \begin{cases}
\sigma^-(M_0) = 0,&\text{if $\tau(M) = M$;}\\
\sigma^-(M_0) - 1,&\text{otherwise.}
\end{cases}
\]

% Now let $\Mbar = (\Mbar,\Fcal,\Vcal)$ any Dieudonn\'e space. We set
% \[
% \tgbar(\Mbar) = \set{m \in \Mbar}{\Vcal(m) \in \Fcal(\Mbar)}.
% \]
% Then we have $\Fcal\Mbar \subset \tgbar(\Mbar)$ and $\Vcal(\tgbar(\Mbar)) \subset \Ker(\Vcal) \subset \tgbar(\Mbar)$. In particular, $\tgbar(\Mbar)$ is a $k[\Fcal,\Vcal]$-submodule of $\Mbar$. Repeating this process we obtain a chain of $k[\Fcal,\Vcal]$-submodules
% \[
% \dots \subset \tgbar^{d+1}(\Mbar) \subset \tgbar^d(\Mbar) \subset \dots \tgbar^0(\Mbar) = \Mbar
% \]
% and we set
% \[
% \tgbar^{\infty}(\Mbar) := \bigcap_{d \geq 0}(\tgbar^d(\Mbar)).
% \]

\begin{proclamation}{Definition}
The integer $\sigma(M) := \sigma^+(M_0)$ is called the {\em superspecial gap of $M$}.
\end{proclamation}

Then $0 \leq \sigma(M) \leq \frac{n-1}{2}$ and if $\sigma(M) > 0$, we have
\begin{equation}\label{GapProp}
[\Lambda^+(M_0) : \Lambda^-(M_0)] = [\Lambda^+(M_1) : \Lambda^-(M_1)] = 2\sigma(M) + 1.
\end{equation}
\end{segment}

\begin{segment}{Ekedahl-Oort strata and the superspecial gap}{EOgap}
The following theorem shows that the (scheme theoretical) Ekedahl-Oort stratification and the (set theoretical) decomposition by the superspecial gap coincide.

\begin{proclamation}{Theorem}\label{EOGap}
Let $x \in \Ncal(k)$, where $k$ is an algebraically closed extension of $\FF_{p^2}$, and let $M$ be the associated unitary Dieudonn\'e module. Let $\sigma$ be an integer with $0 \leq \sigma \leq \frac{n-1}{2}$. Then $x \in \Ncal(\sigma)(k)$ if and only if $\sigma(M) = \sigma$.
\end{proclamation}

\begin{proof}
We have to show that $\sigma(M) = \sigma$ is equivalent to $\Mbar \cong (\MM_{\sigma})_k$ (with the notations of~\eqref{MMSigma}).

For any $k$-vector space $\Ubar$ of $\Mbar_1$ we set
\[
\tgbar(\Ubar) := \Vcal^{-1}(\Fcal(\Ubar)) \cap \Mbar_1,
\]
where $\Vcal^{-1}(\ )$ denotes the preimage under $\Vcal$. We also lift $\tgbar$ by setting
\[
\tau'(U) := (\tau(U) \cap M_1) + pM_1
\]
for any $W(k)$-submodule $U$ of $M_1$ with $pM_1 \subset U$. By~\eqref{SigmaProp} we have $\tau^d(M_1) \supset pM_1$ for all $d \geq 0$ and hence
\[
(\tau')^d(M_1) = \bigcap_{i = 0}^d \tau^i(M_1).
\]

An easy calculation using the explicit definition of $\BB(r)$ and $\SS$ shows that $\Mbar \cong (\MM_{\sigma})_k$ if and only if
\[
\dim(\tgbar^i(\Mbar_1)) = 
\begin{cases}
\dim_k(\Mbar_1) - i,& i \leq \sigma;\\
\dim_k(\Mbar_1) - \sigma,& i \geq \sigma.
\end{cases}
\]

Hence $\Mbar$ and $(\MM_{\sigma})_k$ are isomorphic if and only if $\sigma$ is the smallest integer $d \geq 0$ such that
\[
\bigcap_{i = 0}^d \tau^i(M_1) = \bigcap_{i = 0}^{d+1} \tau^i(M_1).
\]
But by~\eqref{SigmaProp} we have $\sigma(M) = \sigma^-(M_1)$ and therefore $\sigma(M)$ is the smallest integer $d \geq 0$ such that $M_1 \cap \tau(M_1) \cap \dots \cap \tau^d(M_1)$ is $\tau$-invariant. This proves the theorem.
\end{proof}
\end{segment}

%------------------------------------------------------------------

\section{Subschemes of $\Ncal$ indexed by vertices of the \break Bruhat-Tits building of $J$}

In this section we attach to each vertex $\Lambda$ of the Bruhat-Tits building of $J$ a subschemes $\Ncal_{\Lambda}$ of $\Ncal$ and prove that $\Ncal_{\Lambda}$ is isomorphic to a (generalized) Deligne-Lusztig variety.

\begin{segment}{Vertices of the building of $J$}{Building}
We recall now some definitions and results of~\cite{Vo_Uni}. If $\Lambda \subset \Nbf_0$ is a $\ZZ_{p^2}$-lattice we define
\[
\Lambda\vdual = \set{x \in \Nbf_0}{\{x,\Lambda\} \subset \ZZ_{p^2}},
\]
where $\lrcurly$ is the hermitian form on $\Nbf_0$ defined in~\eqref{Curly}. Fix $i \in \ZZ$ such that $ni$ is even and set
\[
\Lcal_i := \set{\text{$\Lambda \subset \Nbf_0$ $\ZZ_{p^2}$-lattice}}{p^{i+1}\Lambda\vdual \subsetneq \Lambda \subset p^i\Lambda\vdual}.
\]
We construct from $\Lcal_i$ an abstract simplicial complex $\Bcal_i$. An $m$-simplex of $\Bcal_i$ is a subset $S \subset \Lcal_i$ of $m+1$ elements which satisfies the following condition. There exists an ordering $\Lambda_0,\dots,\Lambda_m$ of the elements of $S$ such that
\[
p^{i+1}\Lambda\vdual_m \subsetneq \Lambda_0 \subsetneq \Lambda_1 \subsetneq \cdots \subsetneq
\Lambda_m.
\]

Let $\Jtilde = \SU(\Nbf_0,\lrcurly)$ be the derived group of $J$. Note that $\Jtilde$ is a simply connected semisimple algebraic group over $\QQ_p$. There is an obvious action of $\Jtilde(\QQ_p)$ on $\Lcal_i$. We denote by $\Bcal(\Jtilde,\QQ_p)$ the abstract simplicial complex of the Bruhat-Tits building of $\Jtilde$. By \cite{Vo_Uni}~Theorem~3.6 we have a natural identification of $\Bcal_i$ (for fixed $i$) with $\Bcal(\Jtilde,\QQ_p)$. In particular we can identify $\Lcal_i$ with the set of vertices of $\Bcal(\Jtilde,\QQ_p)$. This identification is $\Jtilde(\QQ_p)$-equivariant.

For $\Lambda \in \Lcal_i$ the index of  $p^{i+1}\Lambda\vdual$ in $\Lambda$ is always an odd number, say $2d+1$ by \cite{Vo_Uni}~Remark~2.4. We call the integer $d$ the {\em orbit type of $\Lambda$} and denote it by $t(\Lambda)$. We always have $0 \leq d \leq \frac{n-1}{2}$ and for every such $d$ there exists a $\Lambda \in \Lcal_i$ such that $t(\Lambda) = d$ (loc.~cit.).

Our terminology is explained by the following remark (which follows from \cite{Vo_Uni}~1.17).

\begin{proclamation}{Remark}\label{OrbitType}
Two lattices of $\Lcal_i$ have the same orbit type if and only if they are in the same $\Jtilde(\QQ_p)$-orbit.
\end{proclamation}
\end{segment}

\begin{segment}{$p$-divisible groups attached to a vertex}{VertexPDiv}
Fix again $i \in \ZZ$ such that $ni$ is even. For each $\Lambda \in \Lcal_i$ we will construct two $p$-divisible $O_K$-modules $X_{\Lambda^+}$ and $X_{\Lambda^-}$ over $\FF_{p^2}$ with $O_K$-linear polarization $\lambda_{\Lambda^+}$ and $\lambda_{\Lambda^-}$, respectively. Both will be equipped with an $O_K$-linear quasi-isogeny
\[
\rho_{\Lambda^{\pm}}\colon X_{\Lambda^{\pm}} \to \Xbf
\]
which is compatible with the polarizations up to a power of $p$.

For this we define lattices
\begin{equation}\label{LambdaLattice}
\begin{aligned}
\Lambda_0^+ &:= \Lambda,\\
\Lambda_1^+ &:= \Vbf^{-1}(\Lambda^+_0)\\
\Lambda^+ &:= \Lambda^+_0 \oplus \Lambda^+_1,\\
\Lambda^- &:= \set{x \in \Nbf}{p^{-i}\langle x,\Lambda^+\rangle \subset \ZZ_{p^2}} = p^i(\Lambda^+)^{\perp}.
\end{aligned}
\end{equation}
As $\Fbf = \Vbf$ it follows immediately from the definitions that $\Lambda^{\pm}$ are Dieudonn\'e submodules of the isocrystal $\Nbf$. Clearly $\Lambda^{\pm}_j = \Lambda^{\pm} \cap \Nbf_j$ and therefore the $K$-action on $\Nbf$ restricts to an $O_K$-action on $\Lambda^{\pm}$. As $\Lambda \subset p^i\Lambda\vdual$, the pairing $p^{-i+1}\lrangle$ on $\Nbf$ induces a (not necessarily perfect) $\ZZ_{p^2}$-pairing on $\Lambda^{\pm}$. Hence the unitary $p$-divisible groups $X_{\Lambda^{\pm}}$ associated to the unitary Dieudonn\'e modules $\Lambda^{\pm}$ have all the desired properties. The definition of $\Lambda^+$ shows that $X_{\Lambda^+}$ is a $p$-divisible $O_K$-module of signature $(0,n)$ and an easy calculation shows that this holds for $X_{\Lambda^-}$, too. 

Note that we have
\begin{equation}\label{Lambda0Dual}
\Lambda^-_0 = \set{x \in \Nbf_0}{p^i\langle x, \Vbf^{-1}\Lambda \rangle = p^{-i-1}\{x,\Lambda\} \subset \ZZ_{p^2}} = p^{i+1}(\Lambda^+_0)\vdual
\end{equation}
and hence
\begin{equation}\label{LambdaIndex}
[\Lambda^+_0 : \Lambda^-_0] = [\Lambda^+_1 : \Lambda^-_1] = 2t(\Lambda)+1.
\end{equation}
From this it follows that if $\Lgtilde \in \Lcal_i$ is a second lattice with $\Lgtilde \subsetneq \Lambda$ we have
\begin{equation}\label{TypeIneq}
t(\Lgtilde) < t(\Lambda).
\end{equation}
By definition, $p^{-i}\lrangle$ induces a perfect duality between $\Lambda^+$ and $\Lambda^-$. This duality defines an isomorphism of $p$-divisible $O_K$-modules
$X_{\Lambda^+} \liso {}^tX_{\Lambda^-}$ which makes the diagram
\begin{equation}\label{Dualpm}\xymatrix{
X_{\Lambda^+} \ar[d]_{\rho_{\Lambda^+}} \ar[r]^{\sim} & {}^tX_{\Lambda^-} \\
\Xbf \ar[r]^{\sim}_{\lambda_{\Xbf}} & {}^t\Xbf \ar[u]_{{}^t\rho_{\Lambda^-}}
}\end{equation}
commutative. We set
\begin{equation}\label{BBLambda}
\BB_{\Lambda} := \Lambda^+/\Lambda^-.
\end{equation}
As $p\Lambda^+ \subset \Lambda^-$ and $p^i(\Lambda^+)^{\perp} = \Lambda^-$, $\BB_{\Lambda}$ carries the structure of a unitary Dieudonn\'e space over $\FF_{p^2}$, where the alternating form is induced by $p^{-i+1}\lrangle$.

% \begin{proclamation}{Proposition}
% Let $k$ be an algebraically closed extension of $\FF_{p^2}$. Then $\BB_{\Lambda} \otimes_{\FF_{p^2}} k \cong \BB(2d+1)_k$, where $d$ is the orbit type of $\Lambda$.
% \end{proclamation}
% 
% \begin{proof}
% NOT NEEDED.
% \end{proof}
\end{segment}

\begin{segment}{Schemes $\Ncal_{\Lambda}$ attached to a vertex $\Lambda$}{VertexScheme}
We fix $\Lambda \in \Lcal_i$. For any $\FF_{p^2}$-scheme and for $(X,\rho_X) \in \Ncal_i(S)$ we define quasi-isogenies
\begin{align*}
\rho_{X,\Lambda^+}\colon &X \vartoover{45}{\rho_X} \Xbf_S \vartoover{45}{(\rho_{\Lambda^+})^{-1}_S} (X_{\Lambda^+})_S,\\
\rho_{\Lambda^-,X}\colon &(X_{\Lambda^-})_S \vartoover{45}{(\rho_{\Lambda^-})_S} \Xbf \vartoover{45}{\rho_X^{-1}} X
\end{align*}
Let $t(\Lambda)$ be the orbit type of $\Lambda$. It follows from~\eqref{Dualpm} and~\eqref{LambdaIndex} that
\begin{equation}\label{HtSigma}
\height(\rho_{X,\Lambda^+}) = \height(\rho_{\Lambda^-,X}) = 2t(\Lambda)+1
\end{equation}
and that $\rho_{X,\Lambda^+}$ is an isogeny if and only if $\rho_{\Lambda^-,X}$ is an isogeny.

Let $\Ncal_{\Lambda}$ be the subfunctor of $\Ncal_i \otimes_{\ZZ_{p^2}} \FF_{p^2}$ consisting of those points \mbox{$(X,\rho_X) \in \Ncal_i(S)$} such that $\rho_{X,\Lambda^+}$ or, equivalently, $\rho_{\Lambda^-,X}$ is an isogeny.

\begin{proclamation}{Lemma}
The functor $\Ncal_{\Lambda}$ is representable by a projective \mbox{$\FF_{p^2}$-scheme} and the monomorphism $\Ncal_{\Lambda} \mono \Ncal_i$ is a closed immersion.
\end{proclamation}

\begin{proof}
We set $d := t(\Lambda)$. First of all $\Ncal_{\Lambda} \mono \Ncal_i$ is a closed immersion by \cite{RZ_Period}~2.9. Let $R$ be any $\FF_{p^2}$-algebra. Associating to $(X,\rho_X) \in \Ncal_{\Lambda}(R)$ the kernel of $\rho_{\Lambda^-,X}$ defines a morphism of $\Ncal_{\Lambda}$ to the functor $\Fscr$ on $\FF_{p^2}$-algebras $R$ such that $\Fscr(R)$ is the set of $O_K$-subgroup schemes $\Gcal$ of $X_{\Lambda^-} \otimes R$ of height $2d+1$ such that the polarization $\lambda_{\Lambda^-}$ on $X_{\Lambda^-}$ induces an isomorphism on $(X_{\Lambda^-} \otimes R)/\Gcal$. Conversely, sending $\Gcal \in \Fscr(R)$ to the pair consisting of $X := (X_{\Lambda^-} \otimes R)/\Gcal$ and the inverse $\rho_X$ of quasi-isogeny
\[
\Xbf \otimes R \vartoover{30}{\rho_{\Lambda^-,R}^{-1}} X_{\Lambda^-} \otimes R \lto X
\]
defines an inverse morphism of functors. Therefore $\Ncal_{\Lambda}$ and $\Fscr$ are isomorphic.

The functor $\Fscr$ is a closed subfunctor of the functor parametrizing subgroup schemes of height $2d+1$ of $X_{\Lambda^-}$ which is clearly representable by a projective scheme over $\FF_{p^2}$. Therefore $\Fscr$ and hence $\Ncal_{\Lambda}$ are representable by projective $\FF_{p^2}$-schemes.
\end{proof}

\begin{proclamation}{Lemma}\label{Ncalk}
Let $k$ be an algebraically closed extension of $\FF_{p^2}$. For any lattice $\Gamma$ in $\Nbf$ we set $\Gamma_k = \Gamma \otimes_{\ZZ_{p^2}} W(k)$. For $x \in \Ncal_i(k)$ denote by $M \subset \Nbf_k$ the corresponding unitary Dieudonn\'e module. Then the following assertions are equivalent.
\begin{assertionlist}
\item
$x \in \Ncal_{\Lambda}(k)$.
\item
$M \subset (\Lambda^+)_k$.
\item
$\Lambda^+(M) \subset (\Lambda^+)_k$.
\item
$M_0 \subset \Lambda_k$.
\end{assertionlist}
\end{proclamation}

\begin{proof}
The equivalence of (1) and (2) is clear from the definition of $\Ncal_{\Lambda}$, and the inclusion $M \subset \Lambda^+(M)$ shows that (3) implies (2). As $(\Lambda^+)_k$ is $\tau$-invariant and $\Lambda^+(M)$ is the smallest $\tau$-invariant lattice of $\Nbf_k$ containing $M$ we see that (2) also implies (3). Clearly (2) implies (4) (recall that $\Lambda = \Lambda^+_0$ by definition). Conversely, if $M_0 \subset \Lambda_k$, we have $M_1 \subset p^{-1}\Fcal(M_0) \subset p^{-1}\Fbf(\Lambda)_k = (\Lambda^+_1)_k$.
\end{proof}

This lemma shows that for any algebraically closed extension $k$ of $\FF_{p^2}$ we have
\begin{equation}\label{NcalVcal}
\Ncal_{\Lambda}(k) = \Vcal(\Lambda)(k),
\end{equation}
where $\Vcal(\Lambda)(k)$ is the subset of $\Ncal_i(k)$ defined in \cite{Vo_Uni}~2.3.
\end{segment}

% \begin{intermediate}
% We will now show that $\Ncal_{\Lambda,\red}$ carries a transitive action by the maximal reductive quotient $\Jtilde_{\Lambda}$ of the special fibre of the smooth affine group scheme associated with the vertex of $\Bcal(\Jtilde,\QQ_p)$ which corresponds to $\Lambda$.
% \end{intermediate}

\begin{intermediate}
We will now define a Deligne-Lusztig variety $Y_{\Lambda}$ for a certain reductive group $\Jtilde_{\Lambda}$ over $\FF_p$. Below we will show that $Y_{\Lambda}$ and $\Ncal_{\Lambda}$ are isomorphic. For this we first recall some general facts about Deligne-Lusztig varieties.
\end{intermediate}

\begin{segment}{Deligne-Lusztig varieties}{DLVar}
Let $k$ be a finite field, let $G$ be a connected reductive group over $k$ and let $(W,S)$ be the Weyl system of $G_{\kbar}$, where $\kbar$ is a fixed algebraic closure of $k$. Let $F\colon G \to G$ be the Frobenius morphism with respect to $k$. Then $F$ acts by automorphisms on $W$. As $G$ is quasi-split, $F(S) = S$. For any two subsets $I,I' \subset S$ we denote by ${}^IW^{I'}$ the set of $w \in W$ such that $w$ is the (necessarily unique) element of minimal length in its double coset $W_IwW_{I'}$. Here $W_I$ is the subgroup of $W$ generated by $I$ (and similar for $W_{I'}$).

For any subset $I \subset S$ we denote by $\Par_{G,I} = \Par_I$ the scheme of parabolic subgroups of $G$ of type $I$. It is defined over the extension of degree $r$ of $k$ in $\kbar$, where $r$ is the minimal integer $r \geq 1$ such that $F^r(I) = I$. We denote this field by $\kappa(I)$. Clearly, we have $\kappa(I) = \kappa(F(I))$.

Let $\kappa$ be the compositum of $\kappa(I)$ and $\kappa(I')$ in $\kbar$. The group $G_{\kappa}$ acts diagonally on $\Par_{I,\kappa} \times \Par_{I',\kappa}$. The orbits are in natural bijection to ${}^IW^{I'}$ and we denote by $\Ocal_{I,I'}(w)$ the orbit corresponding to $w \in {}^IW^{I'}$. Clearly, $\Ocal_{I,I'}(w)$ is geometrically irreducible and smooth. Moreover, it is easy to check that
\[
\dim(\Ocal_{I,I'}(w)) = \ell(w) + \dim \Par_{I \cap I'}.
\]
Fix $I \subset S$ and $w \in {}^IW^{F(I)}$. We define the \emph{Deligne-Lusztig variety} $X_I(w)$ as the (locally closed) subscheme of $\Par_{I}$ parametrizing parabolic subgroups $P \in \Par_I$ such that $P$ and $F(P) \in \Par_{F(I)}$ are in relative position $w$.

By definition, $X_I(w)$ is the intersection of $\Ocal_{I,F(I)}(w)$ and the graph of the Frobenius
\[
\Gamma_F\colon \Par_I \to \Par_I \times \Par_{F(I)}, \qquad P \sends (P,F(P)).
\]
It is easily checked that this intersection is transversal. Therefore $X_I(w)$ is smooth of pure dimension
\begin{equation}\label{DimDL}
\dim(X_I(w)) = \ell(w) + \dim \Par_{I \cap F(I)} - \dim(\Par_I).
\end{equation}
Note that $X_I(\id)$ is projective.

Finally, we recall the following theorem of Bonnaf\'e and Rouquier~\cite{BR_IrrDL}.

\begin{proclamation}{Theorem}
The following assertions are equivalent.
\begin{equivlist}
\item
$X_I(w)$ is geometrically irreducible.
\item
$X_I(w)$ is connected.
\item
There exists no $J \subsetneq S$ with $F(J) = J$ and with $W_Iw \subset W_J$.
\end{equivlist}
\end{proclamation}

\begin{proof}
The implication $\text{''(iii) $\Rightarrow$ (i)``}$ is proved in loc.~cit.\ and $\text{''(i) $\Rightarrow$ (ii)``}$ is trivial. Moreover, the proof of Bonnaf\'e and Rouquier (their ''first step``) shows in facts the implication $\text{''(ii) $\Rightarrow$ (iii)``}$.
\end{proof}

Thus $X_I(\id)$ is geometrically irreducible if and only if $\bigcup_{r \geq 1} F^r(I) = S$.
\end{segment}

\begin{segment}{The Deligne-Lusztig variety $Y_{\Lambda}$}{YLambda}
Let $\Lambda \in \Lcal_i$ be a lattice of orbit type $d := t(\Lambda)$. We set $V_{\Lambda} = \Lambda^+_0/\Lambda^-_0$ and we endow $V_{\Lambda}$ with the hermitian pairing $\lrparent = \lrparent_{\Lambda}$ induced by $p^{-i}\lrcurly$. This is a well defined and nondegenerate pairing by~\eqref{Lambda0Dual}. Therefore $(V_{\Lambda},\lrparent)$ is a hermitian space $V_{\Lambda}$ over $\FF_{p^2}$ of dimension $2d+1$.

Set $\Jtilde_{\Lambda} = \SU(V_{\Lambda},\lrparent)$ which is a connected reductive group over $\FF_p$. We remark that $\Jtilde_{\Lambda}$ is the maximal reductive quotient of the special fibre of the smooth affine group scheme attached to the vertex of $\Bcal(\Jtilde,\QQ_p)$ which corresponds to $\Lambda$. We will not use this remark in the sequel.

We denote by $F\colon \Jtilde_{\Lambda} \to \Jtilde_{\Lambda}$ the Frobenius morphism over $\FF_p$ and by $(W,S)$ be the Weyl system of $\Jtilde_{\Lambda}$. Then $F$ induces an automorphism of $W$ which we again denote by $F$. As
\[
\Jtilde_{\Lambda} \otimes_{\FF_p} \FF_{p^2} \cong \SL(V_{\Lambda}) \cong \SL_{2d+1,\FF_{p^2}},
\]
we can identify $(W,S)$ with $(S_{2d+1}, \{s_1,\dots,s_{2d}\})$, where $S_{2d+1}$ is the symmetric group of $\{1,\dots,2d+1\}$ and where $s_i$ is the transposition of $i$ and $i+1$. The automorphism $F$ of $(W,S)$ is induced by the unique nontrivial automorphism of the Dynkin diagram of $\Jtilde_{\Lambda}$, i.e.~$F$ is given by the conjugation with $w_0$, where $w_0$ is the element of maximal length in $S_{2d+1}$ (i.e.~$w_0(i) = 2d+2-i$ for all $i = 1,\dots,2d+1$).

For any subset $I \subset S$, the scheme $\Par_I = \Par_{\Jtilde_{\Lambda},I}$ of parabolic subgroups of $\Jtilde_{\Lambda}$ of type $I$ is defined over $\FF_p$ if and only if $I = F(I)$, otherwise it is defined over $\FF_{p^2}$. We now set
\[
I_{\Lambda} := \{s_1,\dots,s_d,s_{d+2},\dots,s_{2d}\}.
\]
Then $\kappa(I_{\Lambda}) = \FF_{p^2}$. Via the isomorphism $\Jtilde_{\Lambda} \otimes_{\FF_p} \FF_{p^2} \cong \SL(V_{\Lambda})$ we see that $\Par_{I_\Lambda}$ is the $\FF_{p^2}$-scheme which parametrizes $(d+1)$-dimensional subspaces of the $\FF_{p^2}$-vector space $V_{\Lambda}$.

We set
\[
Y_{\Lambda} := X_{I_{\Lambda}}(\id).
\]
By \cite{Vo_Uni}~Lemma~2.17 this definition coincides with the definition given in loc.~cit.~2.12.

\begin{proclamation}{Lemma}\label{PropYLambda}
The $\FF_{p^2}$-scheme $Y_{\Lambda}$ is projective, smooth, and geometrically irreducible of dimension $d = t(\Lambda)$.
\end{proclamation}

\begin{proof}
This is proved in \cite{Vo_Uni}~Proposition~2.18 and Theorem~2.21, but it follows also from general facts on Deligne-Lusztig varieties recalled in~\eqref{DLVar}.
\end{proof}
\end{segment}

\begin{segment}{Modules associated with some isogenies}{LieUnivExt}
To construct an isomorphism $\Ncal_{\Lambda} \iso Y_{\Lambda}$ we will use the Lie algebra of the universal vector extension (see \cite{Me_CrBT}). Let $S$ be an $\FF_p$-scheme. Recall that sending a $p$-divisible group to the Lie algebra of its universal vector extension defines a functor $X \sends D(X)$ from the category of $p$-divisible groups over $S$ to the category of locally free $\Oscr_S$-modules. This functor is compatible with base change $S' \to S$ and we have ${\rm height}(X) = \rk_{\Oscr_S}(D(X))$. Moreover the Hodge filtration is a locally direct summand $H(X)$ of $D(X)$ whose corank is equal to $\rk_{\Oscr_S}(\Lie(X))$.

If $S$ is the spectrum of a perfect field $k$, there is a functorial isomorphism $D(X) \cong M(X)/pM(X) =: \Mbar(X)$, where $M(X)$ is the underlying \mbox{$W(k)$-module} of the covariant Dieudonn\'e module of $X$ and we have \mbox{$H(X) = \Vcal(\Mbar(X))$}.

\begin{proclamation}{Proposition}\label{CokerLie}
Let $\rho\colon X \to Y$ be a homomorphism of $p$-divisible groups over $S$ with $\Ker(\rho) \subset X[p]$. Then $\Coker(D(\rho))$ is a locally free \mbox{$\Oscr_S$-module} and $\rk_{\Oscr_S}(\Coker(D(\rho)) = \height(\rho)$.
\end{proclamation}

Before giving the proof we recall two lemmas.

\begin{proclamation}{Lemma}\label{ComplexExact}
Let $S$ be any scheme and let $\Fscr^1 \to \Fscr^2 \to \Fscr^3$ be a complex of finitely presented $\Oscr_S$-modules such that $\Fscr^2$ and $\Fscr^3$ are locally free. Then $\Fscr\updot$ is exact if and only if for each geometric point $\sbar \to S$ the complex $\Fscr\updot \otimes_{\Oscr_S} \kappa(\sbar)$ is exact.
\end{proclamation}

\begin{proof}
This is a special case of \cite{EGA}~IV,~12.3.3.
\end{proof}

\begin{proclamation}{Lemma}\label{Artin}
Let $R$ be a local Artinian ring with residue field $k$ and let $M$ be a finitely generated $R$-module. Then
\[
\lg_R(M) \leq \lg_R(R)\dim_k(M \otimes_R k).
\]
Moreover, equality holds if and only if $M$ is a free module.
\end{proclamation}

\begin{proof}
This is an immediate application of the Lemma of Nakayama.
\end{proof}

\begin{proof} (of the Proposition)
If $S$ is the spectrum of a perfect field $k$, we have $\dim_k(\Coker(D(\rho))) = \height(\rho)$ by covariant Dieudonn\'e theory. Now for an abritrary $\FF_p$-scheme $S$ the formation of $\Coker(D(\rho))$ is compatible with base change $S' \to S$ and therefore it suffices to show that $\Coker(D(\rho))$ is locally free.

As $\Ker(\rho) \subset X[p]$ there exists a unique isogeny $\rho'\colon Y \to X$ such that $\rho' \circ \rho = p\id_X$ and $\rho \circ \rho' = p\id_Y$. Setting $\varphi := D(\rho)$ and $\varphi' = D(\rho')$ we obtain a complex
\begin{equation}\label{Complex}
\ldots \ltoover{\varphi'} D(X) \ltoover{\varphi} D(Y) \ltoover{\varphi'} D(X) \ltoover{\varphi} \ldots
\end{equation}
whose formation is compatible with base change $S' \to S$. By covariant Dieudonn\'e theory this complex is exact if $S$ is the spectrum of a perfect field. Therefore \eqref{Complex} is exact for arbitrary $S$ by Lemma~\ref{ComplexExact}.

To prove that $\Coker(\varphi)$ is locally free we can assume by standard arguments that $S$ is the spectrum of a local Artin ring $R$. Let $k$ be its residue field. We set
\begin{gather*}
d := \dim_k(\Coker(\varphi \otimes_R \id_k)), \qquad d' := \dim_k(\Coker(\varphi' \otimes_R \id_k)),\\
\ell := \lg_R(R), \qquad h := \rk_R(D(X)) = \rk_R(D(Y)).
\end{gather*}
Then $h = d + d'$ and $h\ell = \lg_R(D(X)) = \lg_R(D(Y))$. We have 
\[
\lg_R(\Coker(\varphi) = \lg_R(\Ker(\varphi)) = \lg_R(\Im(\varphi')) = h\ell - \lg_R(\Coker(\varphi')).
\]
By Lemma~\ref{Artin} we have
\[
h\ell = \lg_R(\Coker(\varphi) + \lg_R(\Coker(\varphi')) \leq d\ell + d'\ell = h\ell
\]
and therefore we must have $\lg_R(\Coker(\varphi) = d\ell$ which implies that $\Coker(\varphi)$ is a free $R$-module again by Lemma~\ref{Artin}.
\end{proof}

\begin{proclamation}{Corollary}\label{KerdS}
Let $S$ be an $\FF_p$-scheme and let $\rho_i\colon X \to Y_i$ for $i = 1,2$ be two isogenies of $p$-divisible groups over $S$ such that $\Ker(\rho_1) \subset \Ker(\rho_2) \subset X[p]$. Then $\Ker(D(\rho_1))$ is locally a direct summand of the locally free $\Oscr_S$-module $\Ker(D(\rho_2))$ and the formation of $\Ker(D(\rho_i))$ commutes with base change $S' \to S$.
\end{proclamation}

\begin{proof}
As $\Coker(D(\rho_i))$ is locally free, $\Im(D(\rho_i))$ is locally a direct summand of $D(Y_i)$ and therefore locally free. Hence the exact sequence
\[
0 \to \Ker(D(\rho_i)) \to D(X) \to \Im(D(\rho_i)) \to 0
\]
locally splits and $\Ker(D(\rho_i))$ is locally a direct summand of $D(X)$ whose formation commutes with base change. Now the exact sequence
\begin{align*}
0 \to \Ker(D(\rho_2))/\Ker(D(\rho_1)) &\lto D(X)/\Ker(D(\rho_1))\\
&\lto D(X)/\Ker(D(\rho_2)) \to 0
\end{align*}
shows that $\Ker(D(\rho_2))/\Ker(D(\rho_1))$ is a locally free $\Oscr_S$-module and therefore $\Ker(D(\rho_1))$ is locally a direct summand of $\Ker(D(\rho_2))$.
\end{proof}
\end{segment}

\begin{segment}{The isomorphism $f\colon \Ncal_{\Lambda} \iso Y_{\Lambda}$}{KeyThm}
% that there is for any $\FF_p$-algebra $R$ a covariant functor $\Dscr$ from the category of $\BT_1$ over $R$ into the category of regular Dieudonn\'e spaces over $R$ (see~\cite{Wd_OS}~(5.5)). If $R$ is a perfect field, $\Dscr$ is isomorphic to the equivalence given by covariant Dieudonn\'e theory.
We will now define a morphism $f\colon \Ncal_{\Lambda} \to Y_{\Lambda}$. Let $R$ be an $\FF_{p^2}$-algebra and let $(X,\rho_X) \in \Ncal_{\Lambda}(R)$. By definition we have isogenies
\[
X_{\Lambda^-,R} \vartoover{35}{\rho_{\Lambda^-,X}} X \vartoover{35}{\rho_{X,\Lambda^+}} X_{\Lambda^+,R}.
\]
The composition is $\rho_{\Lambda} \otimes \id_R$, where $\rho_{\Lambda}\colon X_{\Lambda^-} \to X_{\Lambda^+}$ is the isogeny over $\FF_{p^2}$ which is given via Dieudonn\'e theory by the inclusion $\Lambda^- \subset \Lambda^+$. As $p\Lambda^+ \subset \Lambda^-$, we have $\Ker(\rho_{\Lambda}) \subset X_{\Lambda^-}[p]$ and by~\eqref{BBLambda} the unitary Dieudonn\'e space associated with $\Ker(\rho_{\Lambda})$ is $\BB_{\Lambda}$.

By Corollary~\ref{KerdS}, $E(X) := \Ker(D(\rho_{\Lambda^-,X}))$ is a direct summand of the \mbox{$R$-module} \mbox{$\BB_{\Lambda} \otimes_{\FF_{p^2}} R$}. Its rank is $2d+1$ by~\eqref{HtSigma}. The $O_K$-action on $D$ and on $\BB_{\Lambda}$ defines a decomposition $E(X) = E_0(X) \oplus E_1(X)$ and $\BB_{\Lambda} = \BB_{\Lambda,0} \oplus \BB_{\Lambda,1}$. Note that $\BB_{\Lambda,0} = V_{\Lambda}$ by definition. Therefore $E_0(X)$ is a direct summand of $V_{\Lambda} \otimes_{\FF_{p^2}} R$. We claim that $\rk_R(E_0(X)) = d+1$ and that $E_0(X)^{\perp} \subset E_0(X)$. As $E_0(X)$ is direct summand we can check this after base change to an algebraically closed field and then the claim follows from~\cite{Vo_Uni}~2.12.

We obtain a map $\Ncal_{\Lambda}(R) \to Y_{\Lambda}(R)$ by sending $(X,\rho_X)$ to $E_0(X)$. As the formation of $E(X)$ (and hence of $E_0(X)$) commutes with base change, this map is functorial in $R$ and therefore defines the desired morphism $f$.

\begin{proclamation}{Theorem}\label{Propf}
The morphism $f$ is an isomorphism.
\end{proclamation}

\begin{proof}
By~\cite{Vo_Uni}~2.12 we know that $f$ is a bijection on $k$-valued points for every perfect extension of $\FF_{p^2}$ and hence $f$ is universally bijective. As $\Ncal_{\Lambda}$ is proper and $Y_{\Lambda}$ is separated, $f$ is proper and therefore a universal homeomorphism. In particular we see that $\Ncal_{\Lambda}$ is geometrically irreducible of dimension $d := t(\Lambda)$.

Moreover, we can use Zink's theory of windows for formal $p$-divisible groups \cite{Zi_Wind} and work systematically with a Cohen ring instead of the Witt ring. Then the arguments of loc.~cit.\ show that $f$ is a bijection on $k$-valued points for an arbitrary field extension $k$ of $\FF_{p^2}$.

Now let $k$ be an algebraically closed extension of $\FF_{p^2}$ and let $x \in \Ncal_{\Lambda}(k)$ be a $k$-valued point. We show that the tangent space $T_x(\Ncal_{\Lambda})$ has dimension $\leq d$. This proves the theorem. Indeed, as $\Ncal_{\Lambda}$ is geometrically irreducible of dimension $d$, $\dim(T_x(\Ncal_{\Lambda})) \leq d$ shows that $\Ncal_{\Lambda}$ is smooth. As we already know that $f$ is finite and birational and that $Y_{\Lambda}$ is smooth, Zariski's main theorem shows that $f$ is an isomorphism.

To calculate $T_x(\Ncal_{\Lambda})$ we use Grothendieck-Messing theory. The point $x$ corresponds to a unitary $p$-divisible group $(X,\iota_X,\lambda_X)$ over $k$ endowed with two isogenies
\begin{equation}\label{XX}
X_{\Lambda^-,k} \lto X \lto X_{\Lambda^+,k}.
\end{equation}
Let $\DD^{\pm} = \DD(X_{\Lambda^{\pm},k})$ be the crystal of the Lie algebra of the universal vector extension of $X_{\Lambda^{\pm},k}$ (see \cite{Me_CrBT}~chap.~IV,~\S2). Similarly set $\DD := \DD(X)$. We evaluate this crystals at the PD-thickening $\Spec(k) \mono \Spec(k[\eps]/(\eps^2))$ where the ideal $(\eps)$ is endowed with the trivial PD-structure. We write $\Dtilde^{\pm}$ for the evaluation of $\DD^{\pm}$ and similarly $\Dtilde$ for the evaluation of $\DD$. These three modules are free $k[\eps]/(\eps^2)$-modules of rank $2n$ which an $O_K$-action. Therefore they decompose $\Dtilde^{?} = \Dtilde^{?}_0 \oplus \Dtilde^{?}_1$. Altogether~\eqref{XX} induces homomorphisms of free modules of rank $n$
\begin{equation}\label{DD}
\Dtilde^-_0 \ltoover{\varphi^-} \Dtilde_0 \ltoover{\varphi^+} \Dtilde^+_0.
\end{equation}
As in the construction of $f$ we see that $\Coker(\varphi^-)$ and $\Coker(\varphi^+)$ are free $k[\eps]/(\eps^2)$-modules of rank $d+1$ and rank $d$, respectively. We identify $\Dtilde_0^{\pm}$ with $D(X_{\Lambda^{\pm},k[\eps]/(\eps^2)})_0$ and denote by $\Htilde^{\pm}_0 \subset D(X_{\Lambda^{\pm},k[\eps]/(\eps^2)})_0 = \Dtilde_0^{\pm}$ the zeroth component of the Hodge filtration in $D(X_{\Lambda^{\pm},k[\eps]/(\eps^2)})$. As the signature of $X_{\Lambda^{\pm}}$ is $(0,n)$ by definition~\eqref{VertexPDiv}, we have $\Htilde^{\pm}_0 = \Dtilde^{\pm}_0$. Moreover, let $H_0 \subset D(X)_0$ be the zeroth component of the Hodge filtration given by $X$. As the signature of $X$ equals $(1,n-1)$, $H_0$ is a subspace of codimension $1$.

Now Grothendieck-Messing implies that the tangent space is a subspace of the $k$-vector space of liftings $\Htilde_0$ of $H_0$ to direct summands of $\Dtilde_0$ such that $\varphi^-(\Dtilde^-_0) \subset \Htilde_0$. Hence this tangent space can be identified with a subspace of the tangent space of the projective space parametrizing direct summands of corank $1$ of the module $\Dtilde_0/\Im(\varphi^-)$ which is of rank $d+1$. Therefore we have $\dim(T_x(\Ncal_{\Lambda})) \leq d$.
\end{proof}

\begin{proclamation}{Corollary}\label{PropNLambda}
The closed subscheme $\Ncal_{\Lambda}$ of $\Ncal_i$ is projective, smooth and geometrically irreducible of dimension $t(\Lambda)$.
\end{proclamation}

\begin{proof}
As $Y_{\Lambda}$ has these properties by Lemma~\ref{PropYLambda}, the corollary follows from Theorem~\ref{Propf}.
\end{proof}
\end{segment}

\newpage
%-------------------------------------------------------------------------------------

\section{The global structure of $\Ncal$: The Bruhat-Tits stratification}

\begin{segment}{Combinatorial structure of $\Ncal$}{CombStruc}
We continue to assume that $i$ is an integer such that $ni$ is even, i.e.~that $\Ncal_i \ne \emptyset$.

\begin{proclamation}{Theorem}\label{CombNcalLambda}
Let $\Lambda, \Lgtilde \in \Lcal_i$ be two lattices.
\begin{assertionlist}
\item
Then $\Lambda \subset \Lgtilde$ if and only if $\Ncal_{\Lambda} \subset \Ncal_{\Lgtilde}$. In this case, $t(\Lambda) \leq t(\Lgtilde)$, and $t(\Lambda) = t(\Lgtilde)$ implies $\Lambda = \Lgtilde$.
\item
The following assertions are equivalent.
\begin{equivlist}
\item
$\Lambda \cap \Lgtilde \in \Lcal_i$.
\item
$\Lambda \cap \Lgtilde$ contains a lattice of $\Lcal_i$.
\item
$\Ncal_{\Lambda} \cap \Ncal_{\Lgtilde} \ne \emptyset$.
\end{equivlist}
If these conditions are satisfied,
\[
\Ncal_{\Lambda} \cap \Ncal_{\Lgtilde} = \Ncal_{\Lambda \cap \Lgtilde},
\]
where $\Ncal_{\Lambda} \cap \Ncal_{\Lgtilde}$ denotes the scheme-theoretic intersection in $\Ncal_i$.
\item
The following assertions are equivalent.
\begin{equivlist}
\item
$\Lambda + \Lgtilde \in \Lcal_i$.
\item
$\Lambda + \Lgtilde$ is contained in a lattice of $\Lcal_i$.
\item
$\Ncal_{\Lambda}$ and $\Ncal_{\Lgtilde}$ are both contained in $\Ncal_{\Lambda'}$ for some $\Lambda' \in \Lcal_i$.
\end{equivlist}
If these conditions are satisfied, $\Ncal_{\Lambda + \Lgtilde}$ is the smallest subscheme of the form $\Ncal_{\Lambda'}$ that contains $\Ncal_{\Lambda}$ and $\Ncal_{\Lgtilde}$.
\item
Let $k$ be an algebraically closed extension of $\FF_{p^2}$. Then
\[
\Ncal_i(k) = \bigcup_{\Lambda \in \Lcal_i}\Ncal_{\Lambda}(k).
\]
\end{assertionlist}
\end{proclamation}

\begin{proof}
The nonobvious parts of Assertions (1),~(2), and~(4) follow from \cite{Vo_Uni}~Proposition~2.6. Assertion~(3) is proved as Assertion~(2) by dualizing all lattices.
\end{proof}

\begin{proclamation}{Theorem}\label{NLambdaConn}
Let $i \in \ZZ$ be such that $ni$ is even.
\begin{assertionlist}
\item
The scheme $\Ncal_{i,\red}$ is geometrically connected of pure dimension \mbox{$[\frac{n-1}{2}]$}.
\item
Sending $\Lambda$ to $\Ncal_{\Lambda}$ defines a bijective map
\[
\{\text{$\Lambda \in \Lcal_i$ of orbit type $[\tfrac{n-1}{2}]$}\} \bijective \{\text{irreducible components of $\Ncal_{i,\red}$}\}.
\]
\end{assertionlist}
\end{proclamation}

\begin{proof}
By Theorem~\ref{CombNcalLambda}~(3) the morphism
\begin{equation}\label{NSurj}
\coprod_{\Lambda \in \Lcal_i} \Ncal_{\Lambda} \lto \Ncal_{i,\red} 
\end{equation}
induced by the inclusions $\Ncal_{\Lambda} \mono \Ncal_i$ is surjective.

The lattices $\Lgtilde \in \Lcal_i$ with $\Lambda \subsetneq \Lgtilde$ or $\Lgtilde \subsetneq \Lambda$ are by definition those vertices in the simplicial complex $\Bcal_i$ which are a neighbour of $\Lambda$. Moreover, if $\Lgtilde$ is such a neighbour of $\Lambda$ we have $\Lgtilde \subsetneq \Lambda$ if and only if $t(\Lgtilde) < t(\Lambda)$ by~\eqref{TypeIneq}. In this case we have $\Ncal_{\Lgtilde} \subset \Ncal_{\Lambda}$ by Theorem~\ref{CombNcalLambda}~(1).

As $\Bcal_i$ is isomorphic to $\Bcal(\Jtilde,\QQ_p)$ it follows from the general theory of Bruhat-Tits buildings that for any two vertices $\Lambda$ and $\Lambda'$ of $\Bcal_i$ there exists a sequence of vertices
\[
\Lambda = \Lambda_0, \Lambda_1, \dots, \Lambda_N = \Lambda'
\]
such that $\Lambda_i$ and $\Lambda_{i+1}$ are neighbours and hence either $\Ncal_{\Lambda_i} \subset \Ncal_{\Lambda_{i+1}}$ or $\Ncal_{\Lambda_{i+1}} \subset \Ncal_{\Lambda_i}$. Hence $\Ncal_{\Lambda}$ and $\Ncal_{\Lambda'}$ are contained in the same connected component of $\Ncal_{i,\red}$. As $\Lambda$ and $\Lambda'$ were arbitrary, the surjectivity of~\eqref{NSurj} implies that $\Ncal_{i,\red}$ is connected.

For any $\Lambda \in \Lcal_i$ and any integer $\dtilde$ with $t(\Lambda) \leq \dtilde \leq \frac{n-1}{2}$ there exists a neighbour $\Lgtilde$ of $\Lambda$ with $t(\Lgtilde) = \dtilde$ (\cite{Vo_Uni}~Proposition~2.7). In particular, $\Ncal_{\Lambda}$ is always contained in some $\Ncal_{\Lgtilde}$ with $t(\Lgtilde) = [\frac{n-1}{2}]$. As $\Ncal_{\Lgtilde}$ is irreducible of dimension $t(\Lgtilde)$, this implies that assertion~(2) and also that $\Ncal_{i,\red}$ is of pure dimension $[\frac{n-1}{2}]$.
\end{proof}

For $\Lambda \in \Lcal_i$ we set
\begin{equation}\label{NCal0}
\begin{aligned}
\Lcal_{\Lambda} &:= \set{\Lgtilde \in \Lcal_i}{\Lgtilde \subsetneq \Lambda},\\
\Ncal^0_{\Lambda} &:= \Ncal_{\Lambda} \setminus \bigcup_{\Lgtilde \in \Lcal_{\Lambda}} \Ncal_{\Lgtilde}.
\end{aligned}
\end{equation}

\begin{proclamation}{Proposition}\label{Ncal0Open}
The subset $\Ncal^0_{\Lambda}$ is open and dense in $\Ncal_{\Lambda}$.
\end{proclamation}

\begin{proof}
The set $\Lcal_{\Lambda}$ is a finite set, $\Ncal_{\Lgtilde}$ is a closed subscheme of $\Ncal_{\Lambda}$ for every $\Lgtilde \in \Lcal_{\Lambda}$, and we have
\[
\dim(\Ncal_{\Lgtilde}) = t(\Lgtilde) < t(\Lambda) = \dim(\Ncal_{\Lambda}).
\]
This implies implies the claim.
\end{proof}

Note that by definition we have a disjoint union of locally closed subschemes
\begin{equation}\label{DecompNLambda}
\Ncal_{\Lambda} = \Ncal^0_{\Lambda} \uplus \biguplus_{\Lgtilde \in \Lcal_{\Lambda}} \Ncal^0_{\Lgtilde}.
\end{equation}
We obtain a locally finite stratification $(\Ncal^0_{\Lambda})_{\Lambda \in \Lcal_i}$ of $\Ncal_i$.

\begin{proclamation}{Definition}
The stratification $(\Ncal^0_{\Lambda})_{\Lambda \in \Lcal_i}$ of $\Ncal_i$ is called the \emph{Bruhat-Tits stratification}. The closed subschemes $\Ncal_{\Lambda}$ are called the \emph{closed Bruhat-Tits strata}.
\end{proclamation}

Theorem~\ref{CombNcalLambda} shows that the intersection behaviour of the closed Bruhat-Tits strata can be read off from the simplicial Bruhat-Tits building $\Bcal(\Jtilde,\QQ_p)$. Note that the intersection of two closed Bruhat-Tits strata is always again a closed Bruhat-Tits stratum. In particular it is smooth.

We will now show that the intersection of two $\Ncal_{\Lambda}$'s can be of any dimension within the bounds given by Theorem~\ref{CombNcalLambda}.

\begin{proclamation}{Proposition}\label{InterArb}
Let $d$ and $d'$ be two integers with $0 \leq d, d' \leq (n-1)/2$. Let $\Lambda \in \Lcal_i$ be a lattice with $t(\Lambda) = d$.
\begin{assertionlist}
\item
For any integer $d_-$ with $0 \leq d_- \leq {\rm min}\{d,d'\}$ there exists a lattice $\Lambda' \in \Lcal_i$ with $t(\Lambda') = d'$ such that $\Ncal_{\Lambda} \cap \Ncal_{\Lambda'}$ has dimension $d_-$.
\item
For any integer $d_+$ with ${\rm max}(d,d') \leq d_+ \leq (n-1)/2$ there exists a lattice $\Lambda' \in \Lcal_i$ with $t(\Lambda') = d'$ such that the smallest subscheme $Y$ of the form $\Ncal_{\Lgtilde}$ containing $\Ncal_{\Lambda}$ and $\Ncal_{\Lambda'}$ has dimension $d_+$.
\end{assertionlist}
\end{proclamation}

\begin{proof}
We may assume that $i = 0$. Consider the case that $n$ is even. We can choose a basis $e_1,\dots,e_n$ of $\Nbf_0$ such that $\lrangle$ is given by the matrix $T_{\rm even}$ (see~\eqref{MatrixHerm}). Let $\Lambda(r_1,\dots,r_n)$ be the lattice generated by $p^{r_j}e_j$, $j = 1,\dots,n$, where $r_j \in \ZZ$. Then $\Lambda(r_1,\dots,r_n)\vdual = \Lambda(-r_n-1,-r_{n-1},\dots,-r_1)$. Therefore $\Lambda(r_1,\dots,r_n) \in \Lcal_0$ if and only if
\begin{definitionlist}
\item
$r_1 + r_n = 0$
\item
$0 \leq r_j + r_{n+1-j} \leq 1$ for all $j = 2,\dots,n/2$.
\item
There exists an $j \in \{2,\dots,n/2\}$ such that $r_j + r_{n+1-j} = 0$.
\end{definitionlist}
In this case,
\[
t(\Lambda(r_1,\dots,r_n)) = \#\set{2 \leq j \leq n/2}{r_j + -r_{n+1-j} = 1}.
\]
As $\Jtilde(\QQ_p)$ acts transitively on the set of lattices of type $d$ by Remark~\ref{OrbitType}, we may assume that $\Lambda = \Lambda(0^{n/2},1^d,0^{n/2-d})$. If we define
\[
\Lambda' := \Lambda(0,1^{d'-d_-},0^{n/2-1-d'+d_-},1^{d_-},0^{n/2-d_-})
\]
we see that $\Lambda \cap \Lambda'$ is a lattice in $\Lcal_0$ of type $d_-$.

The case that $n$ is odd and the second assertion can be proved similarly.
\end{proof}

\noindent In particular, given any irreducible component $Y$ of $\Ncal$ and any integer $d$ with $0 \leq d \leq (n-1)/2$, there exists an irreducible component $Y'$ such that $\dim(Y \cap Y') = d$.

\interbreak

\noindent As a second application we examine how many closed Bruhat-Tits strata lie on a given one. For this we define the following invariants. Let $(W, (\,,\,))$ be a nondegenerate $(\FF_{p^2}/\FF_p)$-hermitian space of dimension $l$ (unique up to isomorphism). For any $\FF_{p^2}$-subspace $U$ let $U^\perp$ the orthogonal space with respect to $(\,,\,)$. We fix an integer $r$ with $l/2 \leq r \leq l$ and set
\begin{equation}\label{HermInv}
\begin{aligned}
N(r,W) &:= \set{\text{$U \subset V$ $\FF_{p^2}$-subspace}}{\dim(U) = r, U^\perp \subset U},\\
\nu(r,l) &:= \# N(r,W).
\end{aligned}
\end{equation}
We remark that by \cite{Vo_Uni}~Lemma~2.17, $N(r,W)$ has an interpretation as the set of $\FF_{p^2}$-valued points of a Deligne-Lusztig variety of the group $\SU(W)$.

\begin{remark}{Example}
For all integers $l \geq 2$ we have
\[
\nu(l-1,l) = \#X_{l-1}(\FF_{p^2}),
\]
where $X_l = V_+(X_0^{p+1} + \dots + X_{l}^{p+1})$ is a Fermat hypersurface in $\PP^{l}_{\FF_{p^2}}$. We claim that
\[
\#X_l(\FF_{p^2}) = \begin{cases}
(p^{l+1} + 1)\Sigma_l,&\text{if $l$ is even;}\\
(p^l + 1)\Sigma_l,&\text{if $l$ is odd,}
\end{cases}
\]
where
\[
\Sigma_l := \sum_{j=1}^{[\frac{l-1}{2}]} p^{2j}.
\]
In particular
\[
\nu(1,2) = p+1, \qquad \nu(2,3) = p^3+1.
\]
Let us prove the claim. For this we remark that $\#\mu_{p+1}(\FF_{p^2}) = p+1$, where $\mu_{p+1}$ is the scheme of $(p+1)$-th roots of unity. Therefore, if $a \in \FF^{\times,p+1}_{p^2}$ (i.e.~$a$ is a $(p+1)$-th power in $\FF^{\times}_{p^2}$), there are $p+1$ solutions to the equation $X^{p+1} = a$. Moreover, an element $a \in \FF^{\times}_{p^2}$ lies in $\FF^{\times,p+1}_{p^2}$ if and only if $a \in \FF_p^{\times}$. We set
\[
A_l := \set{x = (x_1,\dots,x_l) \in \AA^l(\FF_{p^2})}{1 + x_1^{p+1} + \dots + x_l^{p+1} = 0}.
\]
If we define $b_l := \#X_l(\FF_{p^2})$ and $a_l := \#A_l$, we have $b_l = b_{l-1} + a_l$ and hence
\begin{equation}\label{Claim1}
b_l = \sum_{j=1}^la_j.
\end{equation}
The set $\set{x \in A_l}{x_l = 0}$ has $a_{l-1}$ elements. The complement is the set $A'_l$ of elements $x$ with $x_l^{p+1} = - 1 - (x_1^{p+1} + \dots + x_{l-1}^{p+1}) \ne 0$. For these we can choose $(x_1,\dots,x_{l-1})$ arbitrary as long as $(x_1,\dots,x_{l-1}) \notin A_{l-1}$ and then we have for $x_l$ still $p+1$ possibilities. Therefore we see that $\#A'_l = (p+1)(p^{2(l-1)} - a_{l-1})$ and hence $a_l = (p+1)p^{2(l-1)} - pa_{l-1}$. An easy induction shows that $a_l = p^{2l-1} + (-1)^{l-1}p^{l-1}$ and then the claim follows from~\eqref{Claim1} by a second easy induction.
\end{remark}

We fix a lattice $\Lambda \in \Lcal_i$ and set $d := t(\Lambda)$. If $\Lambda' \in \Lcal_i$ is a lattice with $\Lambda' \subset \Lambda$ the inclusion induces an $\FF_{p^2}$-linear map $\psi_{\Lambda',\Lambda}\colon V_{\Lambda'} \to V_{\Lambda}$. Now~\cite{Vo_Uni}~2.10 shows that for a fixed $d' \leq d$ we obtain a bijection
\begin{align*}
\set{\Lambda' \in \Lcal_i}{\Lambda' \subset \Lambda, t(\Lambda') = d'} &\leftrightarrow N(d+d'+1,V_{\Lambda}),\\
\Lambda' &\sends \Im(\psi_{\Lambda',\Lambda}).
\end{align*}
Note that $l := \dim V_{\Lambda} = 2d+1$ and hence $l/2 \leq d+d'+1 \leq l$. Therefore Theorem~\ref{CombNcalLambda} shows the following corollary (the second part follows by an easy dualizing argument).

\begin{proclamation}{Corollary}\label{NLLPrime}
\begin{assertionlist}
\item
Let $d'$ be an integer with $0 \leq d' \leq d = t(\Lambda)$. The number of closed Bruhat-Tits strata $\Ncal_{\Lambda'}$ of dimension $d'$ such that $\Ncal_{\Lambda'} \subset \Ncal_{\Lambda}$ is equal to $\nu(d+d'+1, 2d+1)$.
\item
Let $d'$ be an integer with $d = t(\Lambda) \leq d' \leq (n-1)/2$. The number of closed Bruhat-Tits strata $\Ncal_{\Lambda'}$ of dimension $d'$ such that $\Ncal_{\Lambda'} \supset \Ncal_{\Lambda}$ equals $\nu(n - (d+d'+1), n-(2d+1))$.
\end{assertionlist}
\end{proclamation}
\end{segment}

\begin{segment}{Example $n = 4$}{ExGlobal}
As any two non-degenerate hermitian spaces of a fixed dimension over a finite field are isomorphic, $Y_{\Lambda}$ and hence $\Ncal_{\Lambda}$ depends up to isomorphism only on the orbit type $t(\Lambda)$.

\begin{remark}{Example}\label{ExYLambda}
If $t(\Lambda) = 0$, $\Ncal_{\Lambda}$ is a point.

If $t(\Lambda) = 1$, $\Ncal_{\Lambda}$ is isomorphic to the Fermat curve in $\PP^2_{\FF_{p^2}}$ given by the equation $x_0^{p+1} + x_1^{p+1} + x_2^{p+1}$ (\cite{Vo_Uni}~4.11).
\end{remark}

\interbreak

We now apply the general results above to the case $n = 4$. Then for all $i \in \ZZ$, $\Ncal_i$ is nonempty and the $\Ncal_i$ are the connected components of $\Ncal$. $\Ncal_{i, \rm red}$ is equi-dimensional of dimension $1$. The irreducible components are the closed Bruhat-Tits strata $\Ncal_{\Lambda}$ of dimension $1$. Their nonempty intersections are the closed Bruhat-Tits strata $\Ncal_{\Lambda'}$ of dimension $0$.

Each irreducible component $\Ncal_{\Lambda}$ is isomorphic to
\[
V_+(X_0^{p+1} + X_1^{p+1} + X_2^{p+1}) \subset \PP^2_{\FF_{p^2}}.
\]
The intersection points with other irreducible components are precisely the $\FF_{p^2}$-rational points of $\Ncal_{\Lambda}$ and there are $\nu(2,3) = p^3+1$ of them. Through each such intersection point go $p^3+1$ irreducible components and any two of them intersect transversally.
\end{segment}

\begin{segment}{The Ekedahl-Oort strata of $\Ncal_{\Lambda}$}{EONLambda}
Let $i \in \ZZ$ such that $ni$ is even. The Ekedahl-Oort stratification~\eqref{EOStrata} on $\Ncal_i \otimes_{\ZZ_{p^2}} \FF_{p^2}$ induces for every subscheme $S$ a stratification of locally closed subschemes
\[
S = \biguplus_{0 \leq \sigma \leq \frac{n-1}{2}} S(\sigma).
\]

\begin{proclamation}{Theorem}\label{DecompEO}
Fix an integer $\sigma$ with $0 \leq \sigma \leq \frac{n-1}{2}$ and let $\Lambda \in \Lcal_i$. Then
\[
\Ncal_{\Lambda}(\sigma) = \coprod_{\substack{\Lgtilde \in \Lcal_{\Lambda}\\ t(\Lgtilde) = \sigma}} \Ncal_{\Lgtilde}^0.
\]
is a decomposition into open and closed geometrically irreducible subschemes. In particular, $\Ncal_{\Lambda}(\sigma) = \emptyset$ if and only if $\sigma > t(\Lambda)$.
\end{proclamation}

\begin{proof}
By~\eqref{DecompNLambda} it suffices to show the following assertion.
\end{proof}

\begin{proclamation}{Corollary}\label{DecEO}
There exists a decomposition into open and closed geometrically irreducible subschemes
\[
\Ncal_{i,\rm red}(\sigma) = \coprod_{\substack{\Lambda \in \Lcal_i\\ t(\Lambda) = \sigma}} \Ncal^0_{\Lambda}.
\]
\end{proclamation}

\begin{proof}
Let $k$ be an algebraically closed extension of $\FF_{p^2}$. By Theorem~\ref{EOGap} we have
\begin{equation}\label{SetDec}
\Ncal_i(\sigma)(k) = \biguplus_{\substack{\Lambda \in \Lcal_i\\ t(\Lambda) = \sigma}} \Ncal^0_{\Lambda}(k).
\end{equation}
Therefore it suffices to show that $\Ncal^0_{\Lambda}$ is open and closed in $\Ncal_i(\sigma)_{\red}$ for all $\Lambda \in \Lcal_i$ with $t(\Lambda) = \sigma$. By Proposition~\ref{Ncal0Open}, $\Ncal^0_{\Lambda}$ is open in $\Ncal_{\Lambda}$. As $\Ncal_{\Lambda} \cap \Ncal_i(\sigma)_{\red} = \Ncal^0_{\Lambda}$, the subscheme $\Ncal^0_{\Lambda}$ is open in $\Ncal_i(\sigma)_{\red}$. Now it follows from~\eqref{SetDec} that the complement of $\Ncal^0_{\Lambda}$ in $\Ncal_i(\sigma)_{\red}$ is also open and therefore $\Ncal^0_{\Lambda}$ is also closed in $\Ncal_i(\sigma)_{\red}$.
\end{proof}

\begin{proclamation}{Theorem}\label{LocalNcal}
The scheme $\Ncal_{i, \rm red}$ is of pure dimension $[(n-1)/2]$ and locally a complete intersection. The smooth locus of $\Ncal_{i, \rm red}$ is equal to the Ekedahl-Oort stratum $\Ncal_i([(n-1)/2])$.
\end{proclamation}

\begin{proof}
We set $d := [(n-1)/2]$. We already know that all irreducible components are smooth and of dimension $d$. Moreover arbitrary intersections of them are smooth by Theorem~\ref{CombNcalLambda}. Let $Z$ and $Z'$ be two irreducible components with $Z \cap Z' \ne \emptyset$ and set $r := \dim(Z \cap Z')$. There exist irreducible components $Z_{d-1},\dots,Z_{r+1}$ such that $Z \cap Z' \subset Z \cap Z_i$ and $\dim(Z \cap Z_i) = i$ (see Proposition~\ref{InterArb}). This implies that $\Ncal_{i, \rm red}$ is locally a complete intersection and that its singular locus consists of the points which lie on more than a single irreducible component. Therefore its smooth locus is the union of the $\Ncal_{\Lambda}^0$ for all $\Lambda \in \Lcal_i$ with $t(\Lambda) = [(n-1)/2]$ by Theorem~\ref{NLambdaConn}. And this union is $\Ncal_i([(n-1)/2])$ by Corollary~\ref{DecEO}.
\end{proof}

We will now describe the Ekedahl-Oort strata $\Ncal_{\Lambda}(\sigma)$ of $\Ncal_{\Lambda}$ as Deligne-Lusztig varieties for the group $\Jtilde_{\Lambda}$ defined in~\eqref{YLambda}. We set $d = t(\Lambda)$ and identify the Weyl system of $\Jtilde_{\Lambda}$ with $(S_{2d+1}, \{s_1,\dots,s_{2d}\})$ as in~\eqref{YLambda}. As $\Ncal_{\Lambda}(\sigma) = \emptyset$ for $\sigma > d$ by Theorem~\ref{DecompEO}, we assume $\sigma \leq d$. Define
\[
I_{\sigma} := \{s_1,\dots,s_{d-\sigma-1},s_{2d-\sigma+1},\dots,s_{2d}\} \subset S
\]
and $w_{\sigma} \in S_{2d+1}$ as the cycle $(d+\sigma+1, d+\sigma, \dots, d+1)$ which sends $d+\sigma+1$ to $d+\sigma$ etc. Clearly $F(I_{\sigma}) = I_{\sigma}$ and it is easy to check that $w_{\sigma} \in {}^{I_{\sigma}}W^{I_{\sigma}}$. Note that $I_d = I_{d-1} = \emptyset$, but $w_d \ne w_{d-1}$.

\begin{proclamation}{Corollary}\label{EODL}
For all integers $0 \leq \sigma \leq t(\Lambda)$ the isomorphism $\Ncal_{\Lambda} \iso Y_{\Lambda}$ induces an isomorphism of $\Ncal_{\Lambda}(\sigma)$ with the Deligne-Lusztig variety
\[
X_{I_{\sigma}}(w_{\sigma}) \otimes_{\FF_p} \FF_{p^2}.
\]
\end{proclamation}

\begin{proof}
By Theorem~\ref{DecompEO} this is just a reformulation of \cite{Vo_Uni}~Corollary~2.24 (note that in loc.~cit.\ the Deligne-Lusztig variety is defined in a slightly different way and hence the $X_{P_{\sigma}}(w)$ there is the variety $X_{I_{\sigma}}(w^{-1})$ defined here).
\end{proof}
\end{segment}

%--------------------------------------------------------------------------------------

\section{The supersingular locus of the Shimura variety of $\GU(1,n-1)$}\label{Msupersing}

\begin{segment}{The Shimura datum of $\GU(1,n-1)$}{ShimuraDatum}
Let $B$ be a simple $\QQ$-algebra such that $B \otimes_{\QQ} \RR
\cong {\rm M}_m(\CC)$ and such that $B \otimes_{\QQ} \QQ_p \cong {\rm M}_m(K)$ where $K$ is the quadratic unramified extension of $\QQ_p$ fixed in \eqref{PEL}. Note that this implies that $\KK := \Cent(B)$ is a quadratic imaginary extension of $\QQ$ and that $p$ is inert in $\KK$. Let $\star$ be a positive involution on $B$, and let $\VV \ne 0$ be a finitely generated left $B$-module.
We fix a symplectic form $\lrangle\colon \VV \times \VV \to \QQ$ auch that $\langle bv,v'\rangle = \langle v,b\star v'\rangle$ for all $b \in B$ and $v,v' \in \VV$.

We denote by $\GG$ the algebraic group over $\QQ$ of $B$-linear symplectic
similitudes of $(\VV,\lrangle)$. Then $\GG_{\RR}$ is isomorphic to the group of unitary similitudes $\GU(r,s)$ of an hermitian space of signature $(r,s)$ for nonnegative integers $r$ and $s$ with $r+s = n := \dim_{\KK}(\VV)/m$. In particular, $\GG$ is a connected reductive group over $\QQ$. We denote by $E$ the unique subfield of $\CC$ which is isomorphic to $\KK$ and let $\varphi_0$ and $\varphi_1$ be the two isomorphisms $\KK \to E$.

By \cite{Ko_ShFin}~Lemma~4.3 there is a unique $\GG(\RR)$-conjugacy class $X$ of homomorphisms $h\colon \Res_{\CC/\RR}(\GG_{m,\CC}) \to \GG_{\RR}$ such that every $h \in X$ defines a Hodge structure of type $\{(-1,0),(0,-1)\}$ on $\VV$ (with the sign convention of \cite{De_VarSh}) and such that
\[
\VV_{\RR} \times \VV_{\RR} \to \RR, \qquad (v,v') \sends \langle v,h(\sqrt{-1})v'\rangle
\]
is symmetric and positive definite on $\VV_{\RR}$.

Then $(\GG,X)$ is a Shimura datum. Its reflex field is either equal to $\QQ$ (if $r = s$) or equal to $E$ (if $r \ne s$).

We assume that there exist a $\star$-invariant $\ZZ_{(p)}$-order $O_B$ of $B$ such that $O_B \otimes \ZZ_p$ is a maximal order of $B_{\QQ_p}$ and an $O_B$-invariant $\ZZ_p$-lattice $\Gamma$ of $\VV \otimes_{\QQ} \QQ_p$ such that the alternating form on $\Gamma$ induced by $\lrangle$ is a perfect $\ZZ_p$-form. As any maximal order of $M_m(K)$ is conjugate to $M_m(O_K)$, we can and do choose an isomorphism $B \otimes \QQ_p \liso M_m(K)$ which identifies $O_B \otimes \ZZ_p$ with $M_m(O_K)$.
\end{segment}

\begin{segment}{The associated moduli space of abelian varieties}{Shimura}
Denote by $\AA_f^p$ the ring of finite adeles of $\QQ$ with trivial $p$-th component and fix an open compact subgroup $C^p \subset \GG(\AA_f^p)$. Let $\Mcal_{C^p}$ be the moduli space of abelian varieties associated with the data $(B,\star,\VV,\lrangle,O_B,\Gamma,C^p)$ by Kottwitz \cite{Ko_ShFin} \S 5. More precisely, $\Mcal_{C^p}$ is the set-valued functor on the category of $O_{E,(p)}$-schemes which sends each such scheme $S$ to the set of isomorphism classes of tuples $A = (A,\iota_A,\lgbar_A,\hgbar_A)$ where
\begin{bulletlist}
\item
$A$ is an abelian scheme over $S$ of relative dimension equal to $\dim_{\KK}(\VV)$.
\item
$\iota_A\colon O_B \to \End(A) \otimes_{\ZZ} \ZZ_{(p)}$ is a nonzero (and hence injective) homomorphism of $\ZZ_{(p)}$-algebras. Note that if $A\vdual$ is the dual abelian scheme of $A$, we obtain an $O_B$-action $\iota_{A\vdual}\colon O_B \to \End(A\vdual) \otimes_{\ZZ} \ZZ_{(p)}$ on $A\vdual$ by setting
\[
\iota_{A\vdual}(b) := \iota_A(b^*)\vdual.
\] 
\item
$\lgbar_A \subset \Hom(A,A\vdual) \otimes_{\ZZ} \QQ$ is a one-dimensional $\QQ$-subspace which contains a $p$-principal $O_B$-linear polarization.
\item
$\hgbar_A$ is a $C^p$-level structure $\hgbar_A\colon H_1(A,\AA^p_f) \liso V \otimes_{\QQ} \AA^p_f\text{ mod }C^p$.
\end{bulletlist}
such that $(A,\iota_A)$ satisfies Kottwitz's determinant condition of signature $(r,s)$, i.e.~that we have an equality of polynomials
\[
{\rm charpol}(b,\Lie(A)) = (T - \varphi_0(b))^{mr}(T - \varphi_1(b))^{ms} \in \Oscr_S[T]
\]
for all $b \in O_B \cap \KK$ (note that $O_B \cap \KK = O_{\KK,(p)}$ and hence $\varphi_i(O_B \cap \KK) = O_{E,(p)}$; via the structure morphism $S \to \Spec O_{E,(p)}$ we can therefore consider $\varphi_i(b)$ for $b \in O_B \cap \KK$ as sections in $\Oscr_S$).

We call two such tuples $(A,\iota_A,\lgbar_A,\hgbar_A)$ and $(A',\iota_{A'},\lgbar_{A'},\hgbar_{A'})$ isomorphic, if there exists an $O_B$-linear isogeny $\psi\colon A \to A'$ of degree prime to $p$ such that $\psi^*(\lgbar_{A'}) = \lgbar_A$ and $\hgbar_{A'} \circ H_1(\psi,\AA^p_f) = \hgbar_A$.

This functor is represented by a smooth quasi-projective scheme over $O_{E,(p)}$ if $C^p$ is sufficiently small, e.g.~if $C^p$ is contained in a principal congruence subgroup of level $N \geq 3$, where $N$ is an integer prime to $p$. From now on we will assume that this is the case.

We denote by $\Mss_{C^p}$ the supersingular locus of $\Mcal_{C^p} \otimes \FFbar$ considered as a closed reduced subscheme of $\Mcal \otimes \FFbar$, where $\FFbar$ is an algebraic closure of $\FF_{p^2}$. 
\end{segment}

\begin{segment}{Ekedahl-Oort strata of $\Mss_{C^p}$}{EOMss}
We identify $E_p$ with $\QQ_{p^2}$ (and therefore $O_{E_p}$ with $\ZZ_{p^2}$ and the residue field of $O_{E_p}$ with $\FF_{p^2}$). Let $S$ be an $\ZZ_{p^2}$-scheme on which $p$ is locally nilpotent. To every $S$-valued point $A = (A,\iota_A,\lgbar_A,\hgbar_A) \in \Mcal(S)$ we attach the isomorphism class of a unitary $p$-divisible group of signature $(r,s)$ over $S$ (in the sense of~\eqref{UnitaryPDiv}) as follows. Let $X' = A[p^\infty]$ be the $p$-divisible group of $A$. As the ring of endomorphisms of $X'$ is $p$-adically complete, the $O_B$-action on $A$ induces an action of $O_B \otimes_{\ZZ_{(p)}} \ZZ_p = M_m(O_K)$ on $X'$. Morita equivalence tells us that $X' \sends O_K^m \otimes_{M_m(O_K)} X'$ and $X \sends O_K^m \otimes_{O_K} X$ are mutually quasi-inverse functors between the category of $p$-divisible groups $X'$ over $S$ with a left $M_m(O_K)$-action and the category of $p$-divisible groups with $O_K$-action over $S$. We set $X := O_K^m \otimes_{M_m(O_K)} A[p^{\infty}]$ and denote its $O_K$-action by $\iota_X$. Choose a $p$-principal polarization $\lambda_A \in \lgbar_A$. Then $\lambda_A$ induces an $O_K$-linear $p$-principal polarization $\lambda_X$ on $X$. Then $(X,\iota_X,\lambda_X)$ is a unitary $p$-divisible group over $S$ of signature $(r,s)$. Its isomorphism class is independent of the choice of $\lambda_A$.

As in~\eqref{EOStrata} we define the \emph{Ekedahl-Oort strata of $\Mss_{C^p}$}: For this we assume that $r = 1$. Fix an integer $0 \leq \sigma \leq (n-1)/2$. Let $(\Acal,\iota_{\Acal},\lgbar_{\Acal},\hgbar_{\Acal})$ be the restriction of the universal family over $\Mcal_{C^p}$ to $\Mss_{C^p}$. Let $\Xcal$ be the attached isomorphism class of a unitary $p$-divisible group over $\Mss_{C^p}$. Then we denote by $\Mss_{C^p}(\sigma)$ the locally closed subscheme of $\Mss_{C^p}$ such that a morphism $f\colon S \to \Mss_{C^p}$ of $\FF$-schemes factorizes through $\Mss_{C^p}(\sigma)$ if and only if $f^*\Xcal[p]$ is fppf-locally isomorphic to $(\XXbar_{\sigma})_S$. 

The same arguments as in the proof of Theorem~\ref{EOSmooth} show that $\Mss_{C^p}(\sigma)$ is smooth over $\FF$ for all $\sigma$.
\end{segment}

\begin{segment}{The Rapoport-Zink space and the supersingular locus of the Shimura variety}{NcalMcal}
Fix a point $\Abf' = (\Abf',\iota',\lgbar',\hgbar') \in \Mss_{C^p}(\FFbar)$ and let $\Xbf'$ be the attached unitary $p$-divisible group over $\FFbar$ as in~\eqref{EOMss}. As in \eqref{Modp} we can define a functor $\Ncal'$ on the category of $\FFbar$-schemes $S$ where $\Ncal'(S)$ consists of isomorphism classes of pairs $(X,\rho_X)$ where $X$ is a unitary $p$-divisible group over $S$ of signature $(r,s)$ and where $\rho_X\colon X \to \Xbf'_S$ is a quasi-isogeny compatible with the additional structures as explained in~\eqref{Modp}.

We now relate $\Ncal$ with $\Ncal'$ and $\Ncal'$ with $\Mss_{C^p}$. For this we show the following.

\begin{proclamation}{Lemma}\label{SSIsog}
For any two supersingular unitary Dieudonn\'e modules $M$ and $M'$ of signature $(r,s)$ over $\FFbar$ there exists an isogeny $M \to M'$ of unitary Dieudonn\'e modules.
\end{proclamation}

\begin{proof}
Let $H$ be any connected unramified reductive group over $\QQ_p$. We denote by $B(H)$ the set of $\sigma$-conjugacy classes of $H(L)$ where $L$ is the field of fractions of $W(\FFbar)$. Then $B(H)$ classifies isocrystals with $H$-structure, we refer to~\cite{Ko_IsoI} for details.

We consider the hermitian space $V$ introduced in~\eqref{PEL} as a $\QQ_p$-vector space of dimension $2n$ (where $n = r+s$). Let $H := \GL_{\QQ_p}(V)$. The standard representation $G \mono H$ induces a map $B(G) \to B(H)$ and it suffices to show that this map is injective. Let $T$ be a split maximal torus of $H_L$ and let $S$ be a split maximal torus of $G_L$ such that $S \subset T$. We denote by $X_*(S)$ the abelian group of cocharacters of $S$ and set $\Nscr(G) := (X_*(S) \otimes_{\ZZ} \QQ)/W_S$ where $W_S$ is the Weyl group of $(G,S)$. In the same way we define $\Nscr(H)$. Note that $\Nscr(G)$ and $\Nscr(H)$ do not depend (up to unique isomorphism) of the choices of the maximal tori. Kottwitz has defined in~\cite{Ko_IsoI} a map $\nu_G\colon B(G) \to \Nscr(G)$, the Newton map. For $H$ this is simply the map which associates with an isocrystal of height $\dim(V)$ its Newton polygon. We obtain a commutative diagram
\[\xymatrix{
B(G) \ar[r] \ar[d]_{\nu_G} & B(H) \ar[d]^{\nu_H} \\
\Nscr(G) \ar[r] & \Nscr(H).
}\]
As the derived groups of $G$ and $H$ are simply connected, the Newton maps $\nu_G$ and $\nu_H$ are injective and it suffices to show that $\Nscr(G) \to \Nscr(H)$ is injective.

Now $V \otimes_{\QQ_p} L$ is a $K \otimes_{\QQ_p} L$-module. As in~\eqref{Decomp01} we have $K \otimes_{\QQ_p} L = L \times L$ and therefore $V_L = V_0 \times V_1$, and the hermitian pairing on $V_L$ restrict to perfect pairing $V_0 \times V_1 \to L$ which we use to identify $V_1$ with the dual space $V_0^{*}$. Therefore we can identify $G(L)$ with
\[
\set{g = (g_0,g_1) \in \GL(V_0) \times \GL(V_0^{*})}{\exists\ c(g) \in L^{\times}: g_0^*g_1 = c(g)\id_{V_0^*}}.
\]
Hence we have
\begin{align*}
\Nscr(H) &= \set{\lambda = (\lambda_1,\dots,\lambda_{2n}) \in \QQ^{2n}}{\lambda_1 \geq \dots \geq \lambda_n},\\
\Nscr(G) &= \set{\lambda \in \Nscr(H)}{\exists\ \gamma \in \QQ : \lambda_i + \lambda_{2n+1-i} = \gamma\ \forall\ i = 1,\dots,n}
\end{align*}
and $\Nscr(G) \to \Nscr(H)$ is simply the inclusion.
\end{proof}

By the lemma there exists a quasi-isogeny $\tau\colon \Xbf_{\FFbar} \to \Xbf'$ of unitary \mbox{$p$-di}\-visible groups. For any $\FFbar$-scheme $S$ we obtain a bijection
\[
\Ncal(S) \to \Ncal'(S),\qquad (X,\rho_X) \sends (X,\tau \circ \rho_X)
\]
which is functorial in $S$. This defines an isomorphism of formal schemes
\begin{equation}\label{NNPrime}
\zeta\colon \Ncal \otimes_{\ZZ_{p^2}} \FFbar \liso \Ncal'.
\end{equation}

The quasi-isogeny $\tau$ defines an isomorphism of $J$ with the group of auto-quasi-isogenies $J'$ of $\Xbf'$. If we identify $J$ with $J'$ via this isomorphism, the isomorphism $\zeta$ is $J(\QQ_p)$-equivariant.

\interbreak

Now $\Ncal'$ and $\Mss_{C^p}$ are related by the results of~\cite{RZ_Period} as explained in \S~6~of~\cite{Vo_Uni}. We recall the main points (note that the situation considered in~\cite{Vo_Uni} is more special as only the case $n = 3$ and a more special Shimura datum is considered, but the arguments are verbatim the same in this more general case).

Let $I$ be the group of $O_B$-linear of quasi-isogenies in $\End_{O_B}(\Abf') \otimes \QQ$ which respect the space of polarizations $\lgbar'$. This is a reductive group over $\QQ$, which is an inner form of $\GG$. We have isomorphisms of algebraic groups $I_{\QQ_p} \iso J$ and $I_{\AA^p_f} \iso \GG_{\AA^p_f}$ where the second one is given by some $\eta' \in \hgbar'$. In this way we consider $I(\QQ)$ as a subgroup of $J(\QQ_p)$ and of $\GG(\AA_f^p)$.

Denote by $g_1,\dots,g_m \in \GG(\AA^p_f)$ representatives of the finitely many double cosets in $I(\QQ)\backslash \GG(\AA^p_f)/C^p$ and set
\[
\Gamma_j = I(\QQ) \cap g_jC^pg^{-1}_j
\]
for all $j = 1,\dots,m$. Then the subgroup $\Gamma_j \subset J(\QQ_p)$ is discrete and cocompact modulo center. By the assumption made on $C^p$ above, $\Gamma_j$ is torsion free for all $j$ (use Serre's lemma that there is no nontrivial automorphism of an abelian variety which fixes a polarization and the $N$-division points for $N \geq 3$, $(N,p) = 1$).

The uniformization theorem of Rapoport and Zink now provides us with isomorphisms of schemes over $\FFbar$
\[
\coprod_{j=1}^m \Gamma_j\backslash \Ncal'_{\rm red} \cong I(\QQ)\backslash (\Ncal'_{\rm red} \times \GG(\AA^p_f)/C^p) \liso \Mss.
\]
Composition with~\eqref{NNPrime} yields an isomorphism
\begin{equation}\label{NM}
\alpha\colon I(\QQ)\backslash (\Ncal_{\FF,\rm red} \times \GG(\AA^p_f)/C^p) \liso \Mss.
\end{equation} 
Consider the induced surjective morphism
\begin{equation}\label{Psi}
\Psi\colon \coprod_{j=1}^m \Ncal_{\FFbar,\red} \liso \coprod_{j=1}^m \Ncal'_{\rm red} \lto \Mss.
\end{equation}
As explained in the proof of Theorem~6.5 of \cite{Vo_Uni}, our assumption on $C^p$ (and the consequence that $\Gamma_j$ is torsion free) implies that the canonical projection
\[
\pi_j\colon \Ncal'_{\rm red} \to \Gamma_j\backslash \Ncal'_{\rm red}
\]
is a local isomorphism and hence $\Psi$ is a local isomorphism. As $\Gamma_j$ is discrete, the restriction of $\pi_j$ to any closed quasi-compact subscheme of $\Ncal'_{\rm red}$ is finite and therefore the same holds for $\Psi$.
\end{segment}

\begin{segment}{The structure of $\Mss_{C^p}$}{StrucMss}
We now again assume that $r = 1$.

\begin{proclamation}{Theorem}\label{NMss}
Let $C^p$ be as above. The supersingular locus $\Mss_{C^p}$ is of pure dimension $[(n-1)/2]$ and locally of complete intersection. Its smooth locus is the open Ekedahl-Oort stratum $\Mss_{C^p}([(n-1)/2])$.
\end{proclamation}

\begin{proof}
As the morphism $\Psi$ in~\eqref{Psi} is a local isomorphism, the claim follows from Theorem~\ref{LocalNcal} (note that $\Psi$ preserves Ekedahl-Oort strata).
\end{proof}

\noindent Fix an integer $0 \leq d \leq (n-1)/2$. Let $i \in \ZZ$ with $ni$ even. By~\eqref{Building} the group $\Jtilde(\QQ_p)$ acts transitively on the set of lattices in $\Lcal_i$ of orbit type $d$ and therefore on the set of closed Bruhat-Tits strata $\Ncal_{\Lambda} \subset \Ncal_i$ of dimension $d$. Hence Proposition~\ref{JTransConn} implies that $J(\QQ_p)$ acts transitively on the set of all closed Bruhat-Tits strata in $\Ncal$ of dimension $d$. Choose some $\Lambda^{(d)}$ in, say, $\Lcal_0$ such that $t(\Lambda^{(d)}) = d$ and let $C_p^{(d)}$ be its stabilizer in $J(\QQ_p)$. Then we obtain a bijection
\begin{equation}\label{TypedBT}
J(\QQ_p)/C_p^{(d)} \bijective \{\text{closed Bruhat-Tits strata in $\Ncal$ of dimension $d$}\}.
\end{equation}
By Corollary~\ref{PropNLambda} the closed Bruhat-Tits strata $\Ncal_{\Lambda}$ of $\Ncal$ are smooth, geometrically irreducible, and projective. For some $j = 1,\dots,m$, we consider $\Ncal_{\Lambda} \otimes \FF$ as a subscheme in the $j$-th copy of $\Ncal_{\rm red} \otimes \FF$ in the left hand side of~\eqref{Psi}. The restriction of $\Psi$ to this copy is finite. Its image $\Mcal_{\Lambda,j}$ under $\Psi$ is a (geometrically) irreducible and projective subscheme of $\Mss_{C^p}$ of dimension $t(\Lambda)$. We call the subschemes of $\Mss_{C^p}$ of this form the \emph{closed Bruhat-Tits subschemes of $\Mss_{C^p}$}.

Note that the closed Bruhat-Tits subschemes of $\Mss_{C^p}$ of dimension\break \mbox{$[(n-1)/2]$} are the irreducible components of $\Mss_{C^p}$ and that the closed Bruhat-Tits subschemes of dimension $0$ are the superspecial points of $\Mss_{C^p}$. By~\eqref{EONLambda} the closed Bruhat-Tits subschemes of dimension $d$ are the irredu\-cible components of the closure of the Ekedahl-Oort stratum $\Mss_{C^p}(d)$~\eqref{EOMss}. The bijections~\eqref{NM} and~\eqref{TypedBT} therefore show:

\begin{proclamation}{Proposition}\label{IrredMss}
For all integers $0 \leq d \leq (n-1)/2$ there are bijections of finite sets
\begin{align*}
\left\{\begin{matrix}
\text{irreducible components} \\
\text{of $\overline{\Mss_{C^p}(d)}$}
\end{matrix}\right\}
&\bijective
\left\{\begin{matrix}
\text{closed Bruhat-Tits subschemes} \\
\text{of $\Mss_{C^p}$ of dimension $d$}
\end{matrix}\right\}\\
&\bijective
I(\QQ)\backslash (J(\QQ_p)/C_p^{(d)} \times \GG(\AA^p_f)/C^p).
\end{align*}
\end{proclamation}

\noindent Note that $\overline{\Mss_{C^p}([(n-1)/2])} = \Mss_{C^p}$. In particular we get an expression for the number of irreducible components of $\Mss_{C^p}$.

The stabilizer of the connected component $\Ncal_0$ of $\Ncal$ in $J(\QQ_p)$ is
\[
J^0 := \set{g \in J(\QQ_p)}{v_p(c(g)) = 0},
\]
where $c(g)$ is the multiplier of the unitary similitude $g$ and $v_p$ denotes the $p$-adic valuation. As the $\Ncal_i$ are connected by Theorem~\ref{NLambdaConn}, we obtain a bijection
\begin{equation}\label{NConn}
\ZZ = J(\QQ_p)/J^0 \bijective \{\text{connected components of $\Ncal$}\}.
\end{equation}
and therefore the following proposition (where we use that by Theorem~\ref{NLambdaConn} the connected components of $\Ncal$ are geometrically connected).

\begin{proclamation}{Proposition}\label{ConnMss}
There is a bijection of finite sets
\[
\{\text{connected components of $\Mss_{C^p}$}\} \bijective I(\QQ)\backslash (J(\QQ_p)/J^0 \times \GG(\AA^p_f)/C^p).
\]
\end{proclamation}

\noindent The groups $C_p^{(d)}$ are contained in $J^0$ for all $d$. They are special parahoric subgroups of $J$ corresponding to vertices in the Bruhat-Tits building. Via the Propositions~\ref{ConnMss} and~\ref{IrredMss} the fibres of the canonical surjection
\[
I(\QQ)\backslash (J(\QQ_p)/C_p^{(d)} \times \GG(\AA^p_f)/C^p) \to I(\QQ)\backslash (J(\QQ_p)/J^0 \times \GG(\AA^p_f)/C^p)
\]
are in bijection with the set of irreducible components of the Ekedahl-Oort strata of dimension $d$ within a given connected component of $\Mss_{C^p}$.

It can be checked that the number of elements in each of these fibres becomes arbitrarily large if $p$ goes to infinity. In particular we see that the closed Ekedahl-Oort strata of $\Mss_{C^p}$ are highly reducible for large $p$. 

\interbreak

Recall from~\eqref{YLambda} that we defined for each $\Lambda \in \Lcal_i$ ($i$ an integer with $ni$ even) a Deligne-Lusztig variety $Y_{\Lambda}$ whose isomorphism class depended only on $d = t(\Lambda)$. We set $Y_d := Y_{\Lambda}$ (see Example~\ref{ExYLambda} for the shape of $Y_0$ and $Y_1$).

\begin{proclamation}{Proposition}\label{BTIsom}
Fix an integer $0 \leq d \leq (n-1)/2$. For sufficiently small $C^p$, all closed Bruhat-Tits subschemes of $\Mss_{C^p}$ of dimension $d$ are isomorphic to $Y_d$ and therefore smooth. In particular, all irreducible components of $\Mss_{C^p}$ are isomorphic to $Y_{[\frac{n-1}{2}]}$ (and therefore smooth).
\end{proclamation}

\begin{proof}
For sufficiently small $C^p$, the restriction of $\Psi$ to one (or, equivalently, to every) irreducible component $\Ncal_{\Lgtilde}$ (with $t(\Lgtilde) = [(n-1)/2]$) is an isomorphism onto its image by~\cite{Vo_Uni}~6.5. This implies that the same is true for all closed Bruhat-Tits strata $\Ncal_{\Lambda}$ of $\Ncal$ because all of them are contained in some irreducible component. Therefore $\Psi$ induces an isomorphism of $\Ncal_{\Lambda} \cong Y_d$ (with $d = t(\Lambda)$) with its image.
\end{proof}
\end{segment}

%==================================================================

%\bigskip\bigskip

\end{document}